\documentclass[11pt,epsfig]{article}
\usepackage{amssymb}
\usepackage{amsmath}
\usepackage{amsthm}
\usepackage{epsfig}
\usepackage{longtable}
\usepackage{color}
\usepackage{enumitem}
\usepackage[table]{xcolor}

\renewenvironment{proof}[1][Proof]{\noindent\textit{#1. } }{\hfill$\square$}
\newtheoremstyle{thmm}{1.5ex plus 1ex minus .2ex}{1.5ex plus 1ex minus
.2ex}{\rmfamily}{}{\bfseries}{}{1em}{}
\allowdisplaybreaks \textwidth  6.4in \textheight 9in \topmargin
-0.4in \oddsidemargin 0.1in \evensidemargin 0.0in

\theoremstyle{thmm}

\newtheorem{theorem}{Theorem}[section]
\newtheorem{lemma}{Lemma}[section]

\newcommand{\nn}{\nonumber}

\def \endproof{\vrule height8pt width 5pt depth 0pt}
\def\refe#1{(\ref{#1})}

\def\v{\varepsilon}

\def\n{\nu}

\def\d{\delta}

\def\wt{\widetilde}

\def\R{\mathbb{R}}

\def\d{\,{\rm d}}

\def\u{{\bf u}}
\def\U{{\bf U}}
\def\H{{\bf H}}
\def\v{{\bf v}}

\begin{document}
\title{\bf
  New analysis 
of Galerkin-mixed FEMs for incompressible miscible flow in porous media}
\author{
  ~~ Weiwei Sun ~and~Chengda Wu 
  \footnote{Department of Mathematics, City University of Hong Kong,
    Kowloon, Hong Kong.
    The work of the
    authors was supported in part by a grant
    from the Research Grants
    Council of the Hong Kong Special Administrative
    Region, China  
    (Project No. CityU 102613) {\tt maweiw@cityu.edu.hk}, 
{\tt chengda.wu@my.cityu.edu.hk}.  }}
\maketitle

\begin{abstract}
  Analysis of Galerkin-mixed FEMs for incompressible 
  miscible flow in porous media has been investigated extensively 
  in the last several decades.  
  Of particular interest in practical applications is the lowest-order 
  Galerkin-mixed method, { in which a linear Lagrange FE approximation is used 
    for the concentration and the lowest-order Raviart-Thomas FE approximation is used for 
    the velocity/pressure. The previous works only 
    showed the first-order accuracy of the method in $L^2$-norm 
  in spatial direction,} 
  which however is not optimal and valid only under certain extra 
  restrictions on both time step 
  and spatial mesh. 
  In this paper, we provide new and optimal $L^2$-norm error estimates 
  of Galerkin-mixed FEMs for all three components in   
  a general case. In particular, for the lowest-order Galerkin-mixed FEM, 
  we show unconditionally the second-order { accuracy 
  in $L^2$-norm} 
  for the concentration.  Numerical results for both 
  two and three-dimensional models are presented to confirm our 
  theoretical analysis. 
  More important is that our approach can be extended to the 
  analysis of mixed FEMs for many strongly coupled systems 
  to obtain optimal error estimates for all components.

\end{abstract}
\section{Introduction}
\setcounter{equation}{0}
In many engineering applications, incompressible miscible flow in
porous media can be described by
the following miscible displacement system
\begin{align}
  &\Phi\frac{\partial c}{\partial t}-\nabla\cdot(D(\u)\nabla
  c)+\u \cdot\nabla c { + q^P c = \hat c q^I }, 
  \label{e-1}
  \\[3pt]
  &-\nabla\cdot\frac{k(x)}{\mu(c)}\nabla p=q^I-q^P, \label{e-2}
\end{align}
with the initial and boundary conditions:
\begin{align}
  \label{e-3}
  \begin{array}{ll}
    \u\cdot {\bf n}=0,~~
    D(\u)\nabla c\cdot {\bf n}=0
    &\mbox{for}~~x\in\partial\Omega,~~t\in[0,T],\\[3pt]
    c(x,0)=c_0(x)~~ &\mbox{for}~~x\in\Omega,
  \end{array}
\end{align}
where $\u$ denotes the Darcy velocity of the fluid mixture defined by
\begin{align} 
  &\u=-\frac{k(x)}{\mu(c)}\nabla p,  
  \label{e-4}
\end{align}
$p$ is the pressure of the fluid mixture and $c$ is the concentration.
Numerical solutions of the above system \refe{e-1}-\refe{e-4} play a
key role in these applications. Here, $k(x)$ is the
permeability of the medium, $\mu(c)$ is the concentration-dependent
viscosity, $\Phi$ is the porosity of the medium, $q^I$ and $q^P$ are
the given injection and production sources, $\hat c$ is the
concentration in the injection source, and
$D(\u)=[D_{ij}(\u)]_{d\times d}$ is the velocity-dependent
diffusion-dispersion tensor, which may be given in different forms
(see {\cite{BB}} for details). We assume that the system is
defined in a bounded smooth domain $\Omega$ in $\R^d$ $(d=2,3)$ and
$t\in[0,T]$.

In the last several decades, numerous effort has been devoted to
the development of numerical methods for the system \refe{e-1}-\refe{e-4}.
Numerical simulations have been made extensively
in various engineering areas, such as reservoir
simulations and exploration of underground {water and oil
\cite{DFP,Ewing}}. Two review articles for numerical methods
used in these areas {have been written} by Ewing and Wang
\cite{EWang} and Scovazzi et al.\cite{SWML}.  
{The existence of weak solutions of the system was proved
  by Feng \cite{Feng} for the 2D model and by Chen and Ewing \cite{CE}
  for the 3D problem, while the existence of semi-classical/classical 
  solutions is unknown so far. More detailed discussion can be found 
  in \cite{Feng2}. The existence and uniqueness of weak solutions 
  for some different models of Darcy flow were studied in 
  \cite{AFS, RPA}. 
} 
{In \cite{AFS}, the model includes mobile and immobile 
  species, with possibly discontinuous reaction rates, and with a variable 
  porosity that depends on the concentration of the immobile species. 
  The model in \cite{RPA} includes a component for the unsaturated flow 
  (the Richards equation) and another component for reactive transport, with
  nonlinear reaction terms. Whereas here the reaction terms are linear, the 
  diffusion tensor depends on the fluid velocity and the viscosity depends 
  on the solute concentration.
}
Numerical analysis for the system (\ref{e-1})-(\ref{e-4}) in 
two-dimensional space was first presented by Ewing and 
Wheeler \cite{EW} for a standard Galerkin-Galerkin 
approximation $({\cal C}_h, P_h) \in (V_h^r, \widehat V_h^s)$ to
the concentration and pressure in spatial direction,
where $V_h^r$ denotes $C^0$ Lagrange finite element space of
piecewise polynomials of degree $r$ and
$\widehat V_h^s: = V_h^s/\{ constant \}$.
Later, a Galerkin-mixed method was proposed by Douglas et al. \cite{DEW1} for 
solving the system \refe{e-1}-\refe{e-4}. In the Galerkin-mixed method, 
a standard Lagrange type Galerkin approximation ${\cal C}_h \in V_h^r$
was applied for the concentration equation and a mixed approximation
in the Raviart--Thomas finite element space 
($(P_h, \U_h) \in  S_h^k \times \H_h^k$) was used  
for the pressure equation. 
Due to the nature of the continuity of the velocity and the 
discontinuity of the gradient of the pressure in applications, 
the Galerkin-mixed method is more popular in many areas, 
particularly in industries of underground water and oil. 
The error estimate of the semi-discrete Galerkin-mixed method was 
first presented in \cite{DEW1} and later, in \cite{DEW2} for  
a fully discrete semi-implicit Euler scheme. 
In \cite{DEW2}, the error estimate 
\begin{align} 
  \|c^n-{\cal C}_h^n\|_{L^2}+\|p^n-P_h^n\|_{L^2}+\|\u^n - \U_h^n\|_{L^2} 
  \le C (\tau + h_c^{r+1} + h_p^{k+1}) 
  \,  
  \label{error-k} 
\end{align}
was established for $d=2$ under the time step restriction $\tau = o(h)$ and 
an extra spatial mesh size condition,
\begin{align} 
  h_c^{-1} h_p^{k+1} = o(1) 
  \label{mesh-cond-1}
\end{align}  
where $h_c$ and $h_p$ denote the mesh sizes of FE discretization for 
the concentration and pressure equations, respectively. 
Subsequently, some improvements on time step restriction and spatial  
mesh condition were
presented by several authors {\cite{CCW, Dur, LS2, LS3}}.  
In particular, the error estimate \refe{error-k} was proved in 
{\cite{LS2}} 
for $d=2,3$, $h_p=h_c=h$
and $k \ge 1$ in which no time step restriction was required. 
Based on superconvergence analysis, further improvement on the 
spatial mesh condition \refe{mesh-cond-1} was presented in \cite{CWW},  
while the analysis is valid only for regular rectangular meshes. 

The most commonly-used Galerkin-mixed method in practical computation is
the lowest order one ($r=1, k=0$) \cite{CCW, CWW, DEW2, Dur, EWang, SWML}, 
a linear approximation to the concentration and the lowest order
Raviart--Thomas approximation to the pressure and velocity. 
The lowest order Galerkin-mixed method has been widely used in a variety
of numerical simulations and applications, $e.g.$, see  
{\cite{DEW2,Ewing,SWML}}.
In this case,  the error estimate \refe{error-k} reduces to
\begin{align} 
  \| c^n - {\cal C}_h^n \|_{L^2}  +  \|p^n - P_h^n \|_{L^2} 
  + \|\u^n - \U_h^n \|_{L^2} 
  \le C (\tau + h_p + h_c^2)   
  \label{error-0} 
\end{align} 
and the spatial mesh condition \refe{mesh-cond-1} becomes 
\begin{align} 
  h_c^{-1} h_p = o(1)  
  \label{mesh-cond-2} 
  \, . 
\end{align} 

There are two serious concerns arising from previous analysis 
for the popular lowest-order Galerkin-mixed FEM. 
First, it is noted that the error estimate \refe{error-0}
is not optimal for the concentration in $L^2$-norm.
In previous analysis of the lowest-order Galerkin-mixed method, 
a linear approximation to the concentration only produces
the numerical concentration of the accuracy $O(h)$ in spatial direction, 
while the optimal accuracy of a linear approximation 
is $O(h^2)$ in the traditional sense.
Due to the strong coupling of the system, it was assumed that
the one-order lower accuracy of the numerical pressure/velocity may pollute
the numerical concentration through the diffusion-dispersion tensor
$D(\u)$ and the viscosity  $\mu = \mu(c)$.
Our numerical results show that this assumption is incorrect.
Secondly, based on the above spatial mesh condition \refe{mesh-cond-2}, 
one has to use two types of spatial meshes with a much finer one 
for the pressure/velocity equation. 
The Galerkin-mixed method based on the same mesh 
for both concentration and pressure equations ($h_p=h_c$) may not satisfy 
the spatial mesh condition 
\refe{mesh-cond-2} although this method is more efficient and most
commonly-used in practical computation. 

Moreover, mixed finite element methods have been used 
for solving the system \refe{e-1}-\refe{e-4} by combining with
many different schemes in time direction and different approximations
to the concentration equation, such as characteristics type mixed method  
{\cite{AW,Dur,ERW2,FN,KY,Rus1,WSS}},
finite volume method \cite{AO,KNP}, ELLAM \cite{Wang,WLELQ} 
and SUPG \cite{MLG}.    
However, the non-optimality of error estimates for the concentration 
and a time-step/spatial-mesh condition as mentioned
above arise again in the analysis of these methods. 
{Optimal error estimates mentioned above are usually
  based on strong regularity assumptions on data and solution.
  A related topic is the convergence of numerical schemes under
  low regularity assumptions. The convergence of a semidiscrete
  FEM was proved in \cite{CH} for a linear parabolic equation under
  a weak regularity assumption of the solution in 
  $L^2(0,T; H_0^1(\Omega) \cap H^1(0,T;H^{-1}(\Omega))$.
    An implicit Euler scheme with a mixed-DG approximation in spatial
    direction was proposed in \cite{BJM} for the nonlinear system 
    \refe{e-1}-\refe{e-4}. The convergence of the discrete 
    solution to certain weak solution of the system was proved in
    \cite{BJM} by applying the Aubin-Lions compactness on a 
    nonconforming space. Since only some weak regularity assumptions
    on data were made in \cite{RPA}, their analysis implies the
    existence of weak solution of the nonlinear system under 
    weaker assumptions than those in \cite{CE,Feng}. Later, a
    high-order DG scheme in time direction was studied in \cite{RW}. 
    However, no optimal convergence rate/error estimate was obtained 
  under these weak assumptions. }{On the other hand, 
    degenerate cases were analyzed, $e.g.$, 
    in \cite{AW,CHM}. Usually, the convergence rate for degenerate system is
    one-order lower. In particular, in \cite{AW}, a volume corrected 
    characteristics-mixed method is proposed for a purely transport problem. 
    A lower-order $L^1$-norm error estimate $O(h/\sqrt{\tau} + h + \tau^m)$ is
    obtained, where $m$ is related to the accuracy of the characteristic 
    tracing. Different models of coupled Darcy flow and reactive transport 
    are studied, for example, in \cite{AFS, BK, RPA}. These include more 
    complicated, nonlinear reaction terms, but the system is weakly coupled, 
    since the diffusion-dispersion tensor $D({\u}) = I$ and/or the viscosity 
    is constant. An Euler implicit-mixed finite element scheme is analyzed in 
    \cite{RPA}, in which the lowest order mixed FE approximation is used 
    for both the concentration equation and pressure equation. The optimal 
    first-order accuracy is established for the weakly coupled model.
    Numerical methods and analysis for incompressible and immiscible Darcy 
    flow can be found, $e.g.$, in \cite{CHM, Ewing, RKN}.
  }

  The main purpose of this paper is to establish the optimal error estimate
  of Galerkin-mixed methods for all three components, concentration,  
  pressure and velocity, without the time-step restriction
  and spatial mesh condition. In particular, for the lowest order Galerkin-mixed 
  method $(r=1, k=0)$, we will provide the optimal error estimate
  \begin{align} 
    \| c^n - {\cal C}^n_h \|_{L^2}  +  h (\| p^n - P^n_h \|_{L^2} 
  + \|\u^n - \U_h^n) \|_{L^2} ) 
  \le C (\tau + h^2) 
  \label{error-new} 
\end{align}
for $h=h_p=h_c$ unconditionally. 
The analysis is based on an elliptic quasi-projection.
In terms of the projection and a negative norm estimate of Raviart--Thomas
finite element methods for the pressure equation,
the low order accuracy of the velocity will not pollute the concentration
in our analysis and
the lowest order Galerkin-mixed method provides numerical concentration of
the accuracy $O(h^2)$ in $L^2$-norm.
Also we extend our analysis to the general approximation
($(C_h^n, P_h^n, \U_h^n) \in  V_h^r \times S_h^{r-1} \times \H_h^{r-1}$)
to obtain the optimal error estimate
\begin{align} 
  \|c^n-C_h^n\|_{L^2}+h(\|p^n-P_h^n\|_{L^2}+\|\u^n - \U_h^n \|_{L^2}) 
  \le C (\tau + h^{r+1}  )
  \, . 
  \label{error-kk} 
\end{align}
With the numerical concentration $C_h^n$,
a new numerical velocity of the same order accuracy as $C_h^n$ can be
calculated by resolving the (elliptic) pressure equation in a given time level
with a higher-order approximation.
More important is that such a strong coupling can be found in many
other physical systems, $e.g.$, see {\cite{AS,EM,ZYS}}, 
where a higher-order
approximation was also used for one of the computational components.
Our approach can be extended to finite element
analysis for these strongly coupled systems to obtain optimal error estimates
for all components.

The rest of the paper is organized as follows. In Section 2, we
introduce a linearized Euler scheme with 
Galerkin-mixed approximations in the spatial direction for the system
(\ref{e-1})-(\ref{e-4}) and present our main results. In
Section 3, we introduce a new elliptic quasi-projection and 
establish the boundedness of numerical solutions in terms of 
an error splitting technique proposed in {\cite{LS2}}. 
Then we prove the optimal error estimates of  Galerkin-mixed FEMs 
in $L^2$-norm unconditionally. In Section 4, we establish 
some basic estimates of the quasi-projection which were used in the 
proof of the main theorem.  
Finally, numerical simulations for the system in both two and three dimensional
spaces are provided in Section 5. 
Numerical results confirm our theoretical 
analysis that the methods   
provide the optimal accuracy in $L^2$-norm for all three physical 
components. 

\section{The Galerkin-mixed FEM and the main results}
\setcounter{equation}{0}
\subsection{{Notations and assumptions}}
For any integer $m\geq 0$ and $1\leq p\leq\infty$, let $W^{m,p}$ 
be the usual Sobolev spaces and $H^{m}:=W^{m,2}$. In addition, 
we denote by ${\bf H(div)}$ the space of vector-valued functions
${\vec f}\in [L^2(\Omega)]^d$ such that $\nabla\cdot {\vec
f}\in L^2(\Omega)$ {and ${\vec f}\cdot \n = 0$ on $\partial\Omega$}.
{$L^k([0; T];W^{m,p}(\Omega))$ 
  denotes the space of 
  time-dependent functions valued in $W^{m, p}(\Omega)$, which are $L^k$ 
  integrable $w.r.t.$ time in the sense of Bochner, 
  while $H^{1}([0,T];W^{m,n}(\Omega))$ denotes the space of 
  time-dependent functions valued in $W^{m, n}(\Omega)$, which are $H^1$ 
  integrable $w.r.t.$ time in the sense of Bochner. 
}

Let $\mathcal{T}_h$ be a regular triangular partition of $\Omega$ with
$\Omega = \cup_K \Omega_K$ and the mesh size
$h=\max_{\Omega_K \in \mathcal{T}_h} \{ \mathrm{diam} \, \Omega_K\}$.  
For a given division of $\mathcal{T}_h$, we define the classical Lagrange 
finite element spaces by
\begin{align*}
  &V_h^r = \{ v_h \in C^0(\Omega): v_h |_{K} \in P_r(K),\quad \forall 
  K \in \mathcal{T}_h\},
\end{align*}
and Raviart-Thomas finite element spaces \cite{RT,VThomee} by 
\begin{align*}
  &\H_h^s: = \{\v_h\in {\bf H(div)} : \v_h|_K\in [P_s(K)]^d+{\bf x}P_s(K), 
  \quad\forall K\in\mathcal{T}_h\}\nn \\
  &S_h^s: = \{v_h\in L^2 : v_h|_K\in P_s(K), 
  \quad\forall K\in\mathcal{T}_h\}, \quad 
  \widehat S_h^s: = S_h^s/\{constants\},  
\end{align*}
where $P_r(K)$ is the space of polynomials of degree $r\ge 0$ on $K$.
Moreover, we denote by $I_h$ the commonly used Lagrange nodal interpolation 
operator on $V_h^r$.

{Let $\{ t_n \}_{n=0}^{n}$ be a uniform partition in the time direction with 
  the step size $\tau=T/N$ and we denote
  $$
  p^{n} = p(x,t_n),\quad \u^n = \u(x,t_n), \quad c^{n} = c(x,t_n) \, .
  $$
  For any sequence of functions $\{ f^{n} \}_{n=0}^{n}$, we define
  $$
  D_t f^{n+1}=\frac{f^{n+1}-f^{n}}{\tau}\, .
  $$
}
{Throughout this article, we make use of the following assumptions:

  \begin{enumerate}[label=(A\arabic*)]
    \item The solution to the
      initial-boundary value problem (\ref{e-1})-(\ref{e-4}) exists
      and satisfies
      \begin{align}
        \label{StrongSOlEST}
        &\|p\|_{L^\infty([0,T];W^{2,4}\cap H^{s+1})} 
        +\|\u\|_{L^\infty([0,T];W^{1,4}\cap H^{s+1})}
        +\|\u_t\|_{L^2([0,T];H^{s+1})} 
        +\|c\|_{L^\infty([0,T];W^{2,4}\cap H^{r+1})} \nn \\
        &+\|c_t\|_{L^\infty([0,T];H^2)}+\|c_t\|_{L^4([0,T];W^{1,4}\cap H^{r+1})}
        +\|c_{tt}\|_{L^4([0,T];L^4)} \leq K_1.
      \end{align}
    \item $ {\|\hat c\|_{H^1([0,T];H^1(\Omega))}, \,\| q^I \|_{H^1([0,T];H^1(\Omega))}, \, \| q^P \|_{H^1([0,T];H^1(\Omega))}} \leq K_2 $.
    \item $k\in W^{2,\infty}(\Omega)$; 
      $
      k_0^{-1}\leq k(x)\leq k_0, \forall~ x\in\Omega.
      $   
    \item $\mu\in C^1(\R)$; 
      $\exists~\mu_0 >0,~\mu_0^{-1}\leq \mu(s)\leq \mu_0,~\forall~s\in\R.$
    \item Following {\cite{BB}}, the
      diffusion-dispersion tensor is defined by
      \begin{align} 
        D({\u}) = \Phi \left (d_{mt}(|{\u}|) I + 
          d_{lt}(|{\u}|) {\u} \otimes {\u} 
        \right ) ,
      \end{align} 
      where $d_{mt}(z)>d_m>0$, $d_{lt}(z) > 0$ for $z>0$ and $\u \otimes \u 
      =\u \u^T$ denotes a $d \times d$ matrix. 
      For simplicity, here we further assume that $d_{mt}, d_{lt}\in H^3(\R)$ and
      {$\partial_t \nabla D(\v) \in L^{\infty}(0,T;L^2(\Omega))$ 
      for all smooth function $\bf v$} as required in \cite{Whe}. 
      More discussion on  the regularity of the diffusion-dispersion 
      tensor was given in \cite{CLLS,LS3}, where they get optimal $L^p$ error 
      estimates under the regularity assumption of the commonly-used
      Bear--Scheidegger diffusion-dispersion tensor
      $D(\u)\in W^{1,\infty}(\Omega\times(0,T))$ and a similar assumption
      to the exact solution as in \refe{StrongSOlEST}. 
    \item To keep the well-posedness of the initial-boundary value problem
      (\ref{e-1})-(\ref{e-4}), we require
      \begin{equation}
        \int_\Omega q^I \d x=\int_\Omega q^P \d x.  
      \end{equation}
  \end{enumerate}

}

\subsection{{Schemes and main results}}
{Before proposing the fully discrete numerical scheme, 
  we introduce the weak formulation of the deeply coupled system 
  (\ref{e-1})-(\ref{e-4}). 
  Find $p\in L^2(0,T;L^2(\Omega)/\{constants\}) $, $\u\in L^2(0,T;{\bf H(div)})$ 
  and $c\in H^1(0,T;H^1(\Omega))$,  
  such that for all $\v\in {\bf H(div)}$, $\varphi\in L^{2}(\Omega)$ and 
  $\phi\in H^1(\Omega)$, 
  \begin{align}
    & \biggl(\frac{\mu(c)}{k(x)} \u,\,\v\biggl)
    =-(p ,\, \nabla \cdot \v ),
    \\[3pt]
    & (\nabla\cdot \u ,\, \varphi)
    =(q^I-q^P,\, \varphi),
    \\[3pt]
    & (\phi  \partial_{t} c, \, \phi) + (D(\u)
    \nabla c, \, \nabla \phi ) 
    + ( \u\cdot\nabla c,\,
    \phi) {+ (q^P c, \phi) = (\hat c q^I}, \, \phi)
  \end{align}
  {for $a.e.$ $t\in(0,T]$,} where the initial concentration is given by $c(x,0) = c_0(x)$.
}

With the above notations, a fully discrete Galerkin-mixed finite element scheme 
is to find $P_h^{n}\in \widehat S_h^{s}$, $\U_h^{n}\in \H_h^{s}$ 
and ${\cal C}_h^{n}\in V_h^{r}$, $n=0,1,\cdots,N$, 
such that for all $\v_h\in \H_h^{s}$, $\varphi_h\in S_h^{s}$ and 
$\phi_h\in V_h^{r}$,
\begin{align}
  & \biggl(\frac{\mu({\cal C}_h^{n})}{k(x)} \U_h^{n+1},\,\v_h\biggl)
  =-\Big(P_h^{n+1} ,\, \nabla \cdot \v_h \Big),
  \label{e-FEM-1}\\[3pt]
  & \Big(\nabla\cdot \U_h^{n+1} ,\, \varphi_h\Big)
  =\Big(q^I-q^P,\, \varphi_h\Big),
  \label{e-FEM-2}\\[3pt]
  & \Big(\Phi  D_t{\cal C}_h^{n+1}, \, \phi_h\Big) + \Big(D(
    \U_h^{n+1})
  \nabla {\cal C}_h^{n+1}, \, \nabla \phi_h \Big) \nn\\
  &~~~~~~~~~~~~~~~~~~~~~ + \Big( \U_h^{n+1}\cdot\nabla {\cal C}_h^{n},\,
  \phi_h\Big) + 
  {  \Big( q^P {\cal C}_h^{n+1}, \, \phi_h\Big) }
  = \Big(\hat c q^I, \, \phi_h\Big),
  \label{e-FEM-3}
\end{align}
where the initial data ${\cal C}_h^0 = I_h c_0$. Some slightly 
different schemes were investigated 
by several authors. 
{
  In \cite{CCW,DEW2}, a scheme with two different partitions for 
  \refe{e-FEM-1} and \refe{e-FEM-2}-\refe{e-FEM-3} 
  was investigated, while a smaller mesh size was suggested for the pressure/velocity
  than for the concentration for the lowest-order method. 
  In \cite{RPA}, an extra Lipschitz continuous reaction term was 
  introduced for some applications and the lowest-order mixed FEM 
  was used for both concentration equation and pressure equation. 
  A Crank-Nicolson scheme with the coupled convection term 
  $(\U_h^{n+1}\cdot\nabla {\cal C}_h^{n+1}, \phi_h)$ was proposed
  in \cite{CGW}. Moreover, a fully implicit scheme was studied in 
  \cite{ER, Rus1}, where an extra inner iteration  was required at 
  each time step for solving a system of nonlinear equations. In this 
  paper, we only focus our attention to the scheme 
  \refe{e-FEM-1}-\refe{e-FEM-3}, while the analysis for schemes 
  mentioned above is similar. 
}

In this paper, we denote by $C$ a generic positive constant and
by $\epsilon$ a generic small positive constant, which
are independent of $n$, $h$ and $\tau$, and may depend upon $K_1$ , $K_2$ and 
the physical constants $k_0$ and $\mu_0$.
The following classical Gagliardo-Nirenburg interpolation inequality
will be often used in our proof,
\begin{align} 
  \|\partial ^{j} u\|_{L^p} 
  \le C \|\partial^{m} u\|_{L^k}^{a} \, \|u\|_{L^q}^{1-a}
  + C \|u\|_{L^q},
  \label{gn} 
\end{align}
for $0 \le j < m$ and $\frac{j}{m} \le a \le 1$ with $\frac{1}{p} = 
\frac{j}{d} + a \left( \frac{1}{k} - \frac{m}{d}\right)
+(1-a) \frac{1}{q}$,
except $1 < k < \infty$ and $m-j-\frac{n}{k}$ is a non-negative integer,
in which case the above estimate holds only for $\frac{j}{m} \le a < 1$.

We present our main results in the following theorem.
\medskip
\begin{theorem}\label{ErrestFEMSol}
  {\it Suppose that the initial-boundary value problem
    (\ref{e-1})-(\ref{e-4}) {  
    under the assumptions $(A2)-(A6)$} has a unique solution 
    $(p, \u,c)$ which satisfies (\ref{StrongSOlEST}) {
    with $s = r-1$}.  Then there exist positive constants
    $h_0$ and $\tau_0$ such that when $h<h_0$ and $\tau<\tau_0$, the
    finite element system (\ref{e-FEM-1})-(\ref{e-FEM-3}) admits a
    unique solution $(P_h^{n}, \U^{n}_h, {\cal C}_h^{n}) 
    \in( \widehat S_h^{r-1},\H_h^{r-1},V_h^{r})$, $n=1,\cdots,N$, 
    satisfying 
    \begin{align}
      &\max_{1\leq n\leq N}\|{\cal C}_h^{n} -c^{n}\|_{L^2}
      \leq C_0(\tau+h^{r+1}),\\
      &\max_{1\leq n\leq N}\|P_h^{n} - p^{n}\|_{L^2}
      +\max_{1\leq n\leq N}\|
      \U^{n}_h - \u^{n}\|_{L^2} 
      \leq C_0(\tau+h^{r}), \quad \mbox{ for } r \ge 1, 
    \end{align}
    where $C_0$ is a constant, independent of $h$, $\tau$ and $n$, and may be 
    dependent on $K_1$ , $K_2$, $k_0$ and $\mu_0$. 
  }
\end{theorem}
\bigskip
We will present the proof of Theorem 2.1 in the next two sections.

\section{Analysis}\label{SEction3}
\setcounter{equation}{0}
The key to our analysis is a new elliptic quasi-projection. 
In this section, 
we introduce the projection and prove our main results in Theorem 2.1 
in terms of an error splitting technique proposed in {\cite{LS2}}. 
Correspondingly to the fully discrete system \refe{e-FEM-1}-\refe{e-FEM-3},  
we define the time-discrete solution $(P^{n+1}, \U^{n+1}, {\cal C}^{n+1})$  
by the following elliptic system:
\begin{align}
  &\U^{n+1}=-\frac{k(x)}{\mu({\cal C}^{n})}\nabla P^{n+1},
  \label{TDe-fuel-1}
  \\[3pt]
  &\nabla\cdot \U^{n+1}=q^I-q^P, \label{TDe-fuel-2}\\[3pt]
  &\Phi D_t{\cal C}^{n+1} -\nabla\cdot(D( \U^{n+1})\nabla {\cal
  C}^{n+1}) + \U^{n+1} \cdot\nabla{\cal C}^{n}  
  { + q^P {\cal C}^{n+1} } 
  = \hat c q^I, \label{TDe-fuel-4}
\end{align}
for $x\in\Omega$ and $t\in[0,T]$, with the initial and boundary
conditions
\begin{align}
  \label{TDBC}
  \begin{array}{ll}
    \U^{n+1}\cdot {\bf n}=0,~~
    D(\U^{n+1})\nabla {\cal C}^{n+1}\cdot {\bf n}=0
    &\mbox{for}~~x\in\partial\Omega,~~t\in[0,T],\\[3pt]
    {\cal C}^0(x)=c_0(x)~~ &\mbox{for}~~x\in\Omega,
  \end{array}
\end{align}
The condition $\int_\Omega P^{n+1}dx =0 $ is enforced for the uniqueness of
solution. The fully discrete FE solution $(P_h^{n+1}, \U_h^{n+1},{\cal C}_h^{n+1})$ 
can be viewed as a FE solution of the time-discrete system 
\refe{TDe-fuel-1}-\refe{TDBC}.   

\subsection{Preliminary}
Before to prove our main results, we present some lemmas in this section. 
With the solution $(P^{n}, \U^{n}, {\cal C}^{n})$ to the time-discrete system, 
the error functions can be split into
\begin{align}
  & \|\U^{n}_h- \u^{n}\|_{L^2} \le \| \U^{n} - \U_h^{n} \|_{L^2}  
  + \|\U^{n} - \u^{n}
  \|_{L^2},
  \label{sp-1}
  \\
  & \|P_h^{n}- p^{n}\|_{L^2} \le \| P^{n} - P_h^{n} \|_{L^2}  
  + \| P^{n} - p^{n}
  \|_{L^2},
  \label{sp-2}
  \\
  & \|{\cal C}_h^{n} -c^{n} \|_{L^2}  \le \| {\cal C}^{n} - {\cal C}_h^{n} \|_{L^2} 
  + \| {\cal C}^{n} - c^{n} \|_{L^2}.
  \label{sp-3}
\end{align}
The estimates for the second parts of the above 
splitting and the regularity of the 
solution of the time-discrete system (\ref{TDe-fuel-1})-(\ref{TDBC}) were 
given in Theorem 3.1 of \cite{LS2} under a slightly different assumption. 
We present these results 
in the following lemma and the proof is omitted.

\begin{lemma}
  \label{ErrestTDSol}
  {\it Suppose that the initial-boundary value problem
    (\ref{e-1})-(\ref{e-4}) has a unique solution $(p, \u,c)$
    which satisfies (\ref{StrongSOlEST}). Then there exists  
    a positive constant
    $\tau_0^*$ such that when $\tau<\tau_0^*$, the time-discrete system
    (\ref{TDe-fuel-1})-(\ref{TDBC}) admits a unique solution $(P^{n}, \U^n
    , {\cal C}^n )$, $n=1,\cdots,N$, satisfying 
    \begin{align}
      & \|P^{n}\|_{W^{1,4}} +\| \U^{n}\|_{W^{1,4}}
      +\|{\cal C}^{n}\|_{W^{2,4}}+\|D_t{\cal C}^{n}\|_{L^4}  
      +\|\nabla {\cal C}^{n}\|_{L^\infty} \nn\\
      &
      +\biggl(\sum_{n=1}^{N}\tau\|D_tU^{n}\|_{W^{1,4}}^2\biggl)^{\frac{1}{2}}
      +\biggl(\sum_{n=1}^{N}\tau\|D_tC^{n}\|_{H^2}^2\biggl)^{\frac{1}{2}}
      \leq C_1 , 
      \label{StrongSOlESTTD}
    \end{align}
    and
    \begin{align}\label{errorESTTD}
      &\max_{1\leq n\leq N}\| P^{n} - p^{n} \|_{L^4} 
      +\max_{1\leq n\leq N}\| \U^{n} - \u^{n}\|_{L^4} 
      +\max_{1\leq n\leq N}\| {\cal C}^{n} - c^{n}  \|_{L^4}
      \leq C_1 \tau .
    \end{align}
    where $C_1$ is a constant independent of $h$, $\tau$, $n$ and $C_0$ in Theorem 2.1 and 
    may be dependent on $K_1$, $K_2$, $k_0$ and $\mu_0$.
  }
\end{lemma}
\medskip

Now we introduce our elliptic quasi-projection. 
For any fixed integer $n\geq 1$, we denote
by $(\widetilde P^n_h, \widetilde \U_h^n)$ the mixed projection of
$(P^n, \U^n)$ on
$\widehat S_h^{r-1}\times \H_h^{r-1}$ such that 
\begin{align}
  & \biggl(\frac{\mu({\cal C}^{n})}{k(x)}{\wt\U}_h^{n+1},\,\v_h\biggl)
  =-\Big({\wt P}_h^{n+1} ,\, \nabla \cdot \v_h \Big),
  \\[3pt]
  & \Big(\nabla\cdot ({\wt\U}_h^{n+1}-\U^{n+1}) ,\, \varphi_h\Big) =0, \quad 
  \forall (\varphi_h,\v_h)
  \in S_h^{r-1}\times \H_h^{r-1} \, . 
\end{align}
By the classical mixed FE theory \cite{BS, DR, RT, VThomee} and negative norm 
estimates in \cite{JJ}, we have 
\begin{align}
  &\|\U^n-{\wt\U}_h^n\|_{L^p} + \|P^{n}-{\wt P}_h^n\|_{L^p}  \leq Ch, 
  \quad \mbox{for~all}~~2 \le p \le 4,   
  \label{U-Up-Lp}\\
  &\|\U^n-{\wt\U}_h^n\|_{H^{-1}} + \|P^n-{\wt P}_h^n\|_{H^{-1}}  \leq Ch^2, 
  \label{U-Up-H-1} 
\end{align}
For a given $\U^n$, the quasi-projection 
$\Pi_c^n: H^1(\Omega) \rightarrow
V_h^{r}$ is defined by the elliptic problem,
\begin{align}\label{Nonclassical}
  \Big(D(\U^n)\nabla (\Pi_c^{n}{\cal C}^n-{\cal C}^n), \, \nabla \phi_h \Big) 
  +\Big((D({\wt\U}_h^n)-D(\U^n))\nabla {\Pi_c^{n}\cal C}^n , \, \nabla \phi_h \Big)= 0,  
  \nn \\
  \quad \mbox{for~all}~~\phi_h\in V_h^{r},~~n\geq 1,
\end{align}
with $\int_\Omega(\Pi_c^n{\cal C}^n-{\cal C}^n)d x=0$ and $\Pi_c^0:= I_h$.
Clearly, $\Pi_c^n$ is not a projection since $\Pi_c^n{\cal C}_h^n\neq{\cal C}_h^n$,
and reduces to a classical elliptic projection only when $\U_p^n = \U^n$.
We present some basic estimates of the elliptic quasi-projection 
$\Pi_c^n$ in the following lemma and the proof will be given in section 4. 
\begin{lemma}\label{EPro}
  {\it Under the assumptions of Theorem \ref{ErrestFEMSol},  
    there exists $h_1>0$ such that for any $h\leq h_1$ 
    and $2\le p \le 4$
    \begin{align}
      &\|{\cal C}^n-\Pi_c^n {\cal C}^n\|_{L^2} 
      +h\|\nabla ({\cal C}^n-\Pi_c^n {\cal C}^n)\|_{L^p}
      \leq  C_2h^2,   
      \label{p2-1}  \\
      &\|\Pi_c^n{\cal C}^n\|_{W^{1,\infty}} \leq C_2,
      \label{p2-2}
    \end{align}
    and 
    \begin{align}
      &\biggl(\sum_{n=0}^{N-1}\tau\|D_t ({\cal C}^n
      -\Pi_c^n{\cal C}^n)\|_{L^2}^2\biggl)^{1/2}
      \leq  C_2h^2.
      \label{p2-3}
    \end{align}
    where $C_2$ is a constant independent of $h$, $\tau$, $n$, $C_0$ 
    and may be dependent upon 
    $K_1$, $K_2$, $C_1$, $k_0$ and $\mu_0$.
  } 
\end{lemma}

Following the splitting in \refe{sp-1}-\refe{sp-3},
we first present estimates for the first parts
of the splitting in the following lemma.
\begin{lemma}\label{Boundedness}
  {\it  Under the assumptions of Theorem \ref{ErrestFEMSol}, 
    there exist positive constants $\widehat h_0$
    and $\widehat \tau_0$ such that when $h<\widehat h_0$ and 
    $\tau<\widehat \tau_0$, the
    finite element system (\ref{e-FEM-1})-(\ref{e-FEM-3}) admits a
    unique solution $\{P_h^{n}, \U_h^n, {\cal C}_h^n  \}_{n=1}^n 
    \in ( \widehat S_h^{r-1}, \H_h^{r-1}, V_h^r)$,
    which satisfies
    \begin{align}
      &\|{\cal C}_h^{n} -{\cal C}^{n}\|_{L^2} 
      + h(\|P_h^n - P^n\|_{L^2} +\|\U^n_h - \U^n\|_{L^2} ) 
      \leq  C_3 h^2, \quad \mbox{ for } r \ge 1 
      \label{mainerrbd}
    \end{align}
    where $C_3$ is a constant independent of $h$, $\tau$, $n$, $C_0$ and may be 
    dependent upon $K_1$, $K_2$, $C_1$, $C_2$, $k_0$ and $\mu_0$. 
  }
\end{lemma}

\bigskip
\medskip

\noindent{\it Proof}~~~  
Since at each time step 
\refe{e-FEM-1}-\refe{e-FEM-2} is a standard saddle point system 
and the coefficient matrix of the FE system \refe{e-FEM-3} 
is symmetric positive definite, 
the existence and uniqueness of the numerical solution follow immediately.

Let
$$
\theta_p^n = P_h^n - {\wt P}_h^n,\quad
\theta_\u^n = \U^n_h - {\wt\U}_h^n \quad \mbox{ and } \quad
\theta_c^n = {\cal C}_h^n -\Pi_c^n{\cal C}^n.
$$
By noting the projection error estimates in
\refe{U-Up-Lp} and \refe{p2-1}, we only need to prove the following estimate
\begin{align}
  \|\theta_p^{n}\|_{L^2}+\|\theta_\u^{n}\|_{L^2} 
  +\|\theta_c^{n}\|_{L^2}\leq Ch^2.
  \label{error-1}
\end{align}
Since the solution of the time-discrete system
(\ref{TDe-fuel-1})-(\ref{TDBC}) satisfies
\begin{align}
  & \biggl(\frac{\mu({\cal C}^{n})}{k(x)} \U^{n+1},\,v_h\biggl)
  =-\Big(P^{n+1} ,\, \nabla \cdot \v_h \Big),
  \label{tde-FEM-1}\\[3pt]
  & \Big(\nabla\cdot \U^{n+1} ,\, \varphi_h\Big) =\Big(q^I-q^P,\,
  \varphi_h\Big),
  \label{tde-FEM-2}\\[3pt]
  & \Big(\Phi  D_t{\cal C}^{n+1}, \, \phi_h\Big)
  + \Big(D(\U^{n+1})\nabla {\cal C}^{n+1}, \, \nabla \phi_h \Big) \nonumber\\
  &~~~ + \Big( \U^{n+1}\cdot\nabla {\cal C}^{n},\, \phi_h\Big)+ 
  {  \Big(q^P {\cal C}^{n+1}, \, \phi_h\Big) } 
  = \Big(\hat c q^I, \, \phi_h\Big), 
  \label{tde-FEM-3}
\end{align}
for any $\v_h\in \H_h^{r-1}$, $\varphi_h\in S_h^{r-1}$ 
and $\phi_h\in V_h^{r}$, 
from the finite element system
(\ref{e-FEM-1})-(\ref{e-FEM-3}), we can see that the error functions
$(\theta_p^{n}, \theta_\u^{n}, \theta_c^{n})$ satisfy
following system:
\begin{align}
  & \biggl(\frac{\mu({\cal C}^n_h)}{k(x)} \theta_\u^{n+1} + 
    (\frac{\mu({\cal C}^n_h)}{k(x)} 
  - \frac{\mu({\cal C}^{n})}{k(x)}){\wt\U}_h^{n+1}, \v_h\biggl)  
  =-\Big(\theta_p^{n+1},\, \nabla \cdot \v_h \Big),
  \label{erre-FEM-1}\\[3pt]
  & \Big(\nabla\cdot \theta_\u^{n+1} ,\, \varphi_h\Big) =0,
  \label{erre-FEM-2}\\[3pt]
  &
  \Big(\Phi D_t \theta_c^{n+1}, \, \phi_h\Big) + 
  \Big(D(\U_h^{n+1}) \nabla \theta_c^{n+1}, \, \nabla \phi_h \Big) 
  {+ \Big(q^P \theta_c^{n+1},\,   \phi_h\Big)}
  \nn \\ 
  = &\Big(\Phi D_t ({\cal C}^{n+1}-\Pi_c^{n+1}{\cal C}^{n+1}), \,
  \phi_h\Big) -\Big( \U^{n+1}\cdot
    \nabla ({\cal C}_h^{n}-{\cal C}^{n}),\,
  \phi_h\Big)
  \nn \\
  & -\Big(( \U_h^{n+1}- \U^{n+1})\cdot\nabla {\cal C}_h^{n},
  \phi_h\Big)
  -\Big(({\Pi_c^{n+1} {\cal C}^{n+1}}-{\cal C}^{n+1}) q^P, \, \phi_h\Big) 
  \nn \\
  &+\Big((D({\wt\U}_h^{n+1})-D( \U_h^{n+1}))\nabla\Pi_c^{n+1}{\cal C}^{n+1},
  \, \nabla \phi_h \Big) \nn\\
  :=&J_1(\phi_h)+J_2(\phi_h)+J_3(\phi_h)+J_4(\phi_h)+J_5(\phi_h).
  \label{erre-FEM-3}
\end{align}

Taking $\v_h = \theta_\u^{n+1}$ in \refe{erre-FEM-1} leads to 
\begin{align}
  \|\theta_\u^{n+1}\|_{L^2}\leq C\|\theta^{n}_c\|_{L^2}+Ch^2.
  \label{Up-Uh}
\end{align}
Taking $\phi_h= \theta_c^{n+1}$ 
in (\ref{erre-FEM-3}) and by Lemma \ref{EPro}, we get
\begin{align*}
  |J_1(\theta_c^{n+1})| 
  &\leq C(\|\theta_c^{n+1}\|_{L^2}^2
  +\|D_t ({\cal C}^{n+1}-\Pi_c^{n+1}{\cal C}^{n+1})\|_{L^2}^2), 
  \\
  |J_4(\theta_c^{n+1})|&\leq C\|q^P\|_{L^3}
  {\|{\cal C}^{n+1}-\Pi_c^{n+1}{\cal C}^{n+1}\|_{L^2}}\|\theta_c^{n+1}\|_{L^6} 
  \\
  &\leq {\epsilon\|\nabla\theta_c^{n+1}\|_{H^1}^2
  +C\epsilon^{-1}h^4},
\end{align*} 
and by \refe{Up-Uh}, we have 
\begin{align*} 
  |J_5(\theta_c^{n+1})|
  &\leq C\| \nabla \Pi_c^{n+1} {\cal C}^{n+1} \|_{L^\infty}
  \|\theta_u^{n+1} \|_{L^2}\|\nabla \theta_c^{n+1}\|_{L^2} \\
  & \leq \epsilon \|\nabla \theta_c^{n+1}\|_{L^2}^2 
  +\epsilon^{-1}C\| \theta_c^{n} \|_{L^2}^2.
\end{align*}
Moreover, using integration by part and noting the fact that 
$\nabla \cdot \U^{n+1} = q^I - q^P$ 
and $\U^{n+1} \cdot {\bf n} =0$ on the boundary, 
\begin{align*}
  |J_2(\theta_c^{n+1})| & = |( \U^{n+1}\cdot\nabla (\theta_c^{n}
  +\Pi_c^{n}{\cal C}^{n}-{\cal C}^{n}),\, \theta_c^{n+1})|
  \\
  &=| ( (q^I - q^P) (\theta_c^{n}+\Pi_c^{n}{\cal C}^{n}-{\cal C}^{n}),\,  
  \theta_c^{n+1})\\
  &~~~
  + ( \theta_c^{n}+\Pi_c^{n}{\cal C}^{n}-{\cal C}^{n},\,
  \U^{n+1}\cdot\nabla \theta_c^{n+1})|
  \\
  &\le \epsilon\|\theta_c^{n+1}\|_{H^1}^2 
  +C\epsilon^{-1}\|\theta_c^{n}\|_{L^2}^2
  +C\epsilon^{-1}h^4
  \, .
\end{align*}
Finally, we estimate   
\begin{align*}
  |J_3(\theta_c^{n+1})|=&
  | \big(( \U_h^{n+1}- \U^{n+1})\cdot 
  \nabla ({\cal C}_h^{n}-{\cal C}^{n}),\, \theta_c^{n+1}\big) 
  +\big(( \U_h^{n+1}- \U^{n+1})\cdot\nabla {\cal C}^{n},\,  
  \theta_c^{n+1}\big)\big |.
\end{align*} 
By \refe{U-Up-Lp}-\refe{U-Up-H-1} and \refe{Up-Uh},   
\begin{align*}
  & | \big(( \U^{n+1}- \U_h^{n+1})\cdot 
  \nabla ({\cal C}_h^{n}-{\cal C}^{n}),\, \theta_c^{n+1}\big) | 
  \\
  \le &C (\|\theta_u^{n+1}\|_{L^2} +  \|{\wt\U}_h^{n+1}- \U_h^{n+1}\|_{L^2}  )
  (\|\nabla \theta_c^{n}\|_{L^3}
  +\|\nabla (\Pi_c^{n}{\cal C}^{n}-{\cal C}^{n})\|_{L^3})  
  \|\theta_c^{n+1}\|_{L^6} 
  \\
  \le &C (\|\theta^{n}_c\|_{L^2}+h)
  (h^{-d/6}\|\nabla \theta_c^{n}\|_{L^2}+h)  
  \|\theta_c^{n+1}\|_{L^6} 
  \\
  \le &C(h + h^{-d/6} \|\nabla\theta_c^{n} \|_{L^2})( 
    \|\theta_c^{n}\|_{L^2} \|\theta_c^{n+1}\|_{H^1}   
  + h\|\theta_c^{n+1}\|_{H^1}) \,\\
  \le & (Ch^{-d/6} \|\theta^{n}_c \|_{L^2}+Ch +\epsilon ) 
  (\|\theta_c^{n}\|_{H^1}^2+ \|\theta_c^{n+1}\|_{H^1}^2)  
  +C\epsilon^{-1}h^4,
\end{align*} 
and
\begin{align*} 
  & \big |\big(( \U_h^{n+1}- \U^{n+1})\cdot\nabla {\cal C}^{n},\,  
  \theta_c^{n+1}\big)\big | 
  \\
  \le& |\big(\theta_u^{n+1}\cdot\nabla {\cal C}^{n},\,  
  \theta_c^{n+1}\big)\big |+|\big(({\wt\U}_h^{n+1}- \U^{n+1}) 
    \cdot\nabla {\cal C}^{n},\,  
  \theta_c^{n+1}\big)\big |
  \\
  \le& C\|\theta_u^{n+1}\|_{L^2}\|\theta_c^{n+1}\|_{L^2}
  +C\|{\wt\U}_h^{n+1}- \U^{n+1}\|_{H^{-1}} 
  \|\nabla {\cal C}^{n}\theta_c^{n+1}\|_{H^1}
  \\ 
  \le& \epsilon  
  \|\theta_c^{n+1}\|_{H^1}^2+C\|\theta^{n}_c \|_{L^2}^2
  +C\|\theta_c^{n+1} \|_{L^2}^2+C\epsilon^{-1}h^4.
\end{align*} 
Then we further have the estimate 
\begin{align*}
  |J_3(\theta_c^{n+1})|
  &  \le 
  (Ch^{-d/6} \|\theta^{n}_c \|_{L^2}+Ch +\epsilon ) 
  (\|\theta_c^{n}\|_{H^1}^2+ \|\theta_c^{n+1}\|_{H^1}^2)  
  \nn \\ 
  &  +C\|\theta^{n}_c \|_{L^2}^2+C\|\theta_c^{n+1} \|_{L^2}^2+C\epsilon^{-1}h^4. 
\end{align*}
It follows that
\begin{align}
  &\frac{1}{2} D_t \|\theta_c^{n+1}\|^2_{L^2} +
  \big \| \sqrt{D(\U_h^{n+1})}\nabla\theta_c^{n+1}\big \|^2_{L^2}\nn\\
  \leq&
  C(\epsilon+ h^{-d/6}\|\theta^{n}_c\|_{L^2}+h)
  (\|\nabla\theta_c^{n}\|_{L^2}^2 + \|\nabla \theta_c^{n+1}\|_{L^2}^2 ) 
  + C\epsilon^{-1}(\|\theta_c^{n+1}\|_{L^2}^2+\|\theta^{n}_c\|_{L^2}^2)
  + C\epsilon^{-1}h^4\nn\\
  &~~~+ C\big \|D_t ({\cal C}^{n+1}-\Pi_c^{n+1} {\cal C}^{n+1})
  \big \|_{L^2}^2  \, .
  \label{tem-1}
\end{align}
Now we prove the $\tau$-independent estimate
\begin{align}\label{ind}
  \|\theta_c^{n}\|_{L^2}\leq h
\end{align}
from \refe{tem-1} by mathematical induction.
It is easy to see that $\|\theta_c^{0}\|_{L^2}\leq h$,
when $h<h_2$ for some $h_2>0$. We assume that the
inequality (\ref{ind}) holds for $0\leq n\leq k$.
Then there exists a positive constant $h_3$ such that when $h<h_3$,
\refe{tem-1} reduces to
\begin{align}
  & D_t \|\theta_c^{n+1}\|^2_{L^2}  + \big\|\sqrt{D(\U^{n+1}_c)}
\nabla \theta_c^{n+1}) \big \|^2_{L^2}
\nn\\
\leq& Ch^4
+C(\|\theta_c^{n+1}\|_{L^2}^2+\|\theta_c^{n}\|_{L^2}^2)
+C\big\|D_t ({\cal C}^{n+1}-\Pi_c^{n+1} {\cal C}^{n+1}) \big \|_{L^2}^2   
\nn
\end{align}
for $0\leq n\leq k$.
By Gronwall's inequality and \refe{p2-3},
there exists $\tau_1>0$ such that,
\begin{align}
  \|\theta_c^{n+1}\|_{L^2} \leq Ch^2
  \label{nigew3-1}
\end{align}
for $0\leq n\leq k$ and $\tau<\tau_1$, which further implies that 
\begin{align*}
  \|\theta_c^{k+1}\|_{L^2}<h .
\end{align*}
Taking $\widehat \tau_0\leq\min\{\tau_0^*,\tau_1\}$ and
$\widehat h_0\leq\min\{h_1,h_3\}$, 
the mathematical induction is closed and \refe{ind} holds for $1 \le n \le N$. 
Moreover, 
inequalities \refe{Up-Uh} and (\ref{nigew3-1}) hold for all $0\leq n\leq N-1$.

It remains to estimate $\theta_p$.  In a traditional way, we consider
the equation
\begin{align*}
  &-\nabla\cdot\biggl(\frac{k(x)}{\mu({\cal C}^{n})}\nabla g\biggl)
  =\theta_p^{n+1}
\end{align*}
with the boundary condition $\frac{k(x)}{\mu({\cal C}^{n})}\nabla
g\cdot {\bf n}=0$ on $\partial\Omega$. It is easy to see that  
$$ 
\|g\|_{H^2}\leq C\|\theta_p^{n+1}\|_{L^2}. 
$$  
Let
$$
\v_h=Q_h\biggl(\frac{k(x)}{\mu({\cal C}^{n})}\nabla g\biggl)
$$
where $Q_h: {\bf H(div)} \rightarrow \H_h^{r-1}$ is a
projection such that \cite{VThomee} for ${\bf w} \in {\bf H(div)}$,  
\begin{align}
  \Big(\nabla\cdot({\bf w}-Q_h{\bf w})\, , \chi_h \Big) = 0 ,\quad
  \mbox{for~all}~~\chi_h \in S_h^{r-1} .  
  \label{Q}
\end{align}
Then
$$
(\varphi_h, \nabla \cdot \v_h)  =-
(\varphi_h, \theta_p^{n+1}),~~\mbox{ for } \varphi_h \in S_h^{r-1}
$$
and from (\ref{erre-FEM-1}) and the classical result
$\|Q_h w\|_{L^2}\leq C\|w\|_{H^1}$, we obtain
\begin{align*}
  \|\theta_p^{n+1}\|_{L^2}^2 & = \biggl(\frac{\mu({\cal
    C}^{n}_h)}{k(x)} \U_h^{n+1}  
    -\frac{\mu(\Pi_c^{n}{\cal C}^{n})}{k(x)}{\wt\U}_h^{n+1},
  Q_h\biggl(\frac{k(x)}{\mu({\cal C}^{n})}\nabla g\biggl) \biggl)
  \nn \\
  & \leq C(\|\theta^{n}_c\|_{L^2}+\|\theta_u^{n+1}\|_{L^2})
  \biggl\|\frac{k(x)}{
  \mu({\cal C}^{n})}\nabla g\biggl\|_{H^1}\\
  &\leq Ch^2\|\theta_p^{n+1}\|_{L^2} ,
\end{align*}
which implies that
\begin{align}
  \|\theta_p^{n+1}\|_{L^2} \leq Ch^2 .
  \label{Pp-Ph}
\end{align}
\refe{error-1} follows \refe{Up-Uh}, \refe{nigew3-1} and \refe{Pp-Ph}. 
The proof of Lemma \ref{Boundedness} is completed. 
\endproof

\subsection{Proof of Theorem \ref{ErrestFEMSol}}
For $r=1$, Theorem \ref{ErrestFEMSol} can be proved by combining 
Lemma \ref{ErrestTDSol} and Lemma \ref{Boundedness}.  
In this section, we only prove the case $r \ge 2$.

From Lemma \ref{ErrestTDSol}, Lemma \ref{EPro},  \refe{U-Up-Lp}, \refe{error-1} and 
the Gagliardo-Nirenburg inequality \refe{gn}, 
we can see the boundedness of numerical solution 
\begin{align}
  \|P_h^n\|_{L^{\infty}}+\|\U_h^n\|_{L^{\infty}} 
  +\|{\cal C}_h^n\|_{W^{1,6}}\le C.
  \label{boundedness}
\end{align}

Similarly, for any fixed integer $n\geq 1$ 
we denote by $({\wt p}_h^n, {\wt\u}_h^n)$ the classical mixed 
projection of $(p^n, \u^n)$ on $(\widehat S_h^{r-1}, \H_h^{r-1})$   
such that 
\begin{align}
  & \biggl(\frac{\mu(c^{n+1})}{k(x)}{\wt\u}_h^{n+1},\, \v_h\biggl)
  =-\Big({\wt p}_h^{n+1} ,\, \nabla \cdot \v_h \Big),
  \\[3pt]
  & \Big(\nabla\cdot ({\wt\u}_h^{n+1}-\u^{n+1}) ,\, \varphi_h\Big) =0, \quad 
  \mbox{ for all } (\varphi_h,v_h) \in S_h^{r-1}\times \H_h^{r-1}.
\end{align}
By the classical mixed method theory \cite{BS, DR, VThomee} and negative
norm estimates in \cite{JJ}, we have
\begin{align}
  &\|\u^n-{\wt\u}_h^n\|_{L^2} + \|p^n-{\wt p}_h^n\|_{L^2}  \leq Ch^r, 
  \label{nU-Up-Lp}\\
  &\|\u^n-{\wt \u}_h^n\|_{H^{-1}} + \|p^n-{\wt p}_h^n\|_{H^{-1}}  \leq Ch^{r+1}, 
  \label{nU-Up-H-1} 
\end{align}
and by inverse inequalities and noting $r \ge 2$, 
\begin{align}
  \|{\wt\u}_h^{n}\|_{L^{\infty}}\leq C.
  \label{BoundProU}
\end{align}
For a given $\u^n$, we define 
an elliptic quasi-projection ${\wt \Pi}_c^n: H^1(\Omega) \rightarrow
V_h^{r}$,  slightly different from one in section 3.1, by 
\begin{align}\label{nNonclassical}
  \Big(D(\u^n)\nabla ({\wt \Pi}_c^n c^n-c^{n}), \, \nabla \phi_h \Big) 
  + \Big((D({\wt \u}_h^n)-D(\u^n))\nabla {\wt \Pi}_c^{n}c^n , \, \nabla \phi_h \Big)= 0, \nn \\
  \quad \mbox{for~all}~~\phi_h\in V_h^{r},~~n\geq 1,
\end{align}
with $\int_\Omega({\wt \Pi}_c^{n}c^{n}-c^{n})d x=0$ and ${\wt \Pi}_c^0:= I_h$.
By a proof similar to Lemma \ref{EPro}, we can get basic estimates 
of the elliptic quasi-projection ${\wt \Pi}_c^{n}$ as follows.
\begin{lemma}\label{nEpro}
  {\it Under the assumptions of Theorem \ref{ErrestFEMSol},  
    there exists $\widehat h_1>0$ such that for any $h\leq\widehat h_1$,
    \begin{align}
      &\|c^{n}-{\wt \Pi}_c^{n} c^{n}\|_{L^2} 
      +h\|\nabla (c^{n}-{\wt \Pi}_c^{n} c^{n})\|_{L^2}
      \leq  Ch^{r+1}, 
      \label{np2-1}  \\
      &\|{\wt \Pi}_c^{n}c^{n}\|_{W^{1,\infty}} \leq C,
      \label{np2-2}
    \end{align}
    and 
    \begin{align}
      &\biggl(\sum_{n=0}^{N-1}\tau\|D_t (c^{n}
      -{\wt \Pi}_c^{n}c^{n})\|_{L^2}^2\biggl)^{1/2}
      \leq  Ch^{r+1}.
      \label{np2-3}
    \end{align}
  } 
\end{lemma}
Now we start to prove Theorem \ref{ErrestFEMSol}. 
Let
$$
{\wt \theta}_p^{n} = P_h^{n} - {\wt p}_h^{n},\quad
{\wt \theta}_u^{n} = \U_h^{n} - {\wt \u}_h^{n} \quad \mbox{ and } \quad
{\wt \theta}_c^{n} = {\cal C}_h^{n} -{\wt \Pi}_c^{n}c^{n}.
$$
We prove below the estimate
\begin{align}
  \|{\wt \theta}_p^{n}\|_{L^2}+\|{\wt \theta}_\u^{n}\|_{L^2}+\|{\wt \theta}_c^{n}\|_{L^2} 
  \leq C(\tau + h^{r+1}).
  \label{nerror-1}
\end{align}
From \refe{e-1}-\refe{e-3} and the finite element system \refe{e-FEM-1}-\refe{e-FEM-3},
the error function $({\wt \theta}_p^{n}$, ${\wt \theta}_u^{n}$, ${\wt \theta}_c^{n})$ satisfies
\begin{align}
  & \biggl(\frac{\mu({\cal C}_h^{n})}{k(x)} {\wt \theta}_\u^{n+1} + 
    (\frac{\mu({\cal C}_h^{n})}{k(x)}-\frac{\mu(c^{n+1})}{k(x)}){\wt \u}_h^{n+1}, 
  \v_h\biggl) =-\Big({\wt \theta}_p^{n+1},\,
  \nabla \cdot \v_h \Big),
  \label{nerre-FEM-1}\\[3pt]
  & \Big(\nabla\cdot {\wt \theta}_\u^{n+1} ,\, \varphi_h\Big) =0,
  \label{nerre-FEM-2}\\[3pt]
  &
  \Big(\Phi D_t {\wt \theta}_c^{n+1}, \, \phi_h\Big) + 
  \Big(D(\U_h^{n+1}) \nabla {\wt \theta}_c^{n+1}, \, \nabla \phi_h \Big) 
  {+\Big( q^P{\wt \theta}_c^{n+1}, \, \phi_h\Big)} 
  \nn \\
  = &(T^{n+1}_c, \phi_h)+\Big(\Phi D_t (c^{n+1}-{\wt \Pi}_c^{n+1}c^{n+1}), \,
  \phi_h\Big) -
  \Big( \u^{n+1}\cdot \nabla ({\cal C}_h^{n}-c^{n}),\, \phi_h\Big)
  \nn \\
  & -\Big(( \U_h^{n+1}- \u^{n+1})\cdot\nabla {\cal C}_h^{n},\,\phi_h\Big)
  -\Big({({\wt \Pi}_c^{n+1}c^{n+1}-c^{n+1})} q^P, \, \phi_h\Big) \nn \\
  &+\Big((D({\wt \u}_h^{n+1})-D( \U_h^{n+1}))\nabla{\wt \Pi}_c^{n+1}c^{n+1},
  \, \nabla \phi_h \Big) \nn\\
  :=&(T^{n+1}_c, \phi_h)+{\wt J}_1(\phi_h)+{\wt J}_2(\phi_h)+{\wt J}_3(\phi_h)
  +{\wt J}_4(\phi_h)+{\wt J}_5(\phi_h),
  \label{nerre-FEM-3}
\end{align}
where $T^{n+1}_c$ denotes the truncation error.
By the regularity assumption (\ref{StrongSOlEST}), we have 
\begin{align}
  \sum_{k = 1}^{n}\tau\|T^k_c\|_{L^2}^2\leq C\tau^2.
\end{align}

Letting $v_h = {\wt \theta}_\u^{n+1}$ in \refe{nerre-FEM-1} and 
by \refe{BoundProU}, we further have 
\begin{align}
  \|{\wt \theta}_\u^{n+1}\|_{L^2}\leq C\|{\wt \theta}^{n}_c\|_{L^2} 
  + Ch^{r+1} + C\tau.   
  \label{nUp-Uh}
\end{align}
Taking $\phi_h= {\wt \theta}_c^{n+1}$ in (\ref{nerre-FEM-3}) and using Lemma \ref{nEpro} gives 
\begin{align*}
  |{\wt J}_1({\wt \theta}_c^{n+1})| 
  &\leq C(\|{\wt \theta}_c^{n+1}\|_{L^2}^2
  +\|D_t (c^{n+1}-{\wt \Pi}_c^{n+1}c^{n+1})\|_{L^2}^2), 
  \\
  |{\wt J}_4({\wt \theta}_c^{n+1})|&\leq C\|q^P\|_{L^3}
  {\|c^{n+1}-{\wt \Pi}_c^{n+1}c^{n+1}\|_{L^2}}\|{\wt \theta}_c^{n+1}\|_{L^6} 
  \\
  &\leq {\epsilon\|{\wt \theta}_c^{n+1}\|_{H^1}^2
  +C\epsilon^{-1} h^{2r+2}},
\end{align*} 
and by \refe{nUp-Uh}, we get 
\begin{align*} 
  |{\wt J}_5({\wt \theta}_c^{n+1})|&\leq C\| \nabla {\wt \Pi}_c^{n+1} c^{n+1} \|_{L^\infty}
  \|{\wt \theta}_\u^{n+1} \|_{L^2}\|\nabla {\wt \theta}_c^{n+1}\|_{L^2} \\
  & \leq C\epsilon \|\nabla {\wt \theta}_c^{n+1}\|_{L^2}^2
  +\epsilon^{-1}\| {\wt \theta}_c^{n} \|_{L^2}^2 + C\epsilon^{-1}(h^{2r+2}+\tau^2),
\end{align*}
Moreover,  using integration by part and noting the fact that 
$\nabla \cdot \u^{n+1} = q^I - q^P$ 
and $\u^{n+1} \cdot {\bf n} =0$ on the boundary, 
\begin{align*}
  |{\wt J}_2({\wt \theta}_c^{n+1})| &  
  = |( \u^{n+1}\cdot\nabla ({\wt \theta}_c^{n}
  +{\wt \Pi}_c^{n}c^{n}-c^{n}),\, {\wt \theta}_c^{n+1})|
  \\
  &=| ( (q^I - q^P) ({\wt \theta}_c^{n}+{\wt \Pi}_c^{n}c^{n}-c^{n}),\,  
  {\wt \theta}_c^{n+1})\\
  &~~~
  + ( {\wt \theta}_c^{n}+{\wt \Pi}_c^{n}c^{n}-c^{n},\,
  \u^{n+1}\cdot\nabla {\wt \theta}_c^{n+1})|
  \\
  &\le \epsilon\|{\wt \theta}_c^{n+1}\|_{H^1}^2 
  +C\epsilon^{-1}\|{\wt \theta}_c^{n}\|_{L^2}^2
  +C\epsilon^{-1}h^{2r+2}
  \, .
\end{align*}
Finally, we rewrite $\wt J_3$ by   
\begin{align*}
  |{\wt J}_3({\wt \theta}_c^{n+1})|=&
  | \big(( \U_h^{n+1}- \u^{n+1})\cdot \nabla ({\cal C}_h^{n}-c^{n}),\, 
  {\wt \theta}_c^{n+1}\big) 
  +\big(( \U_h^{n+1}- \u^{n+1})\cdot\nabla c^{n},\,  
  {\wt \theta}_c^{n+1}\big)\big |.
\end{align*} 
By \refe{nU-Up-Lp}-\refe{nU-Up-H-1} and \refe{nUp-Uh},   we have 
\begin{align*}
  & | \big(( \u^{n+1}- \U_h^{n+1})\cdot 
  \nabla ({\cal C}_h^{n}-c^{n}),\, {\wt \theta}_c^{n+1}\big) | 
  \\
  \le &C (\|\theta_\u^{n+1}\|_{L^3} +  \|{\wt\U}_h^{n+1}- \U^{n+1}\|_{L^3} 
  + \|\U^{n+1}- \u^{n+1}\|_{L^3}  )
  \nn \\
  &(\|\nabla {\wt \theta}_c^{n}\|_{L^2}
  +\|\nabla ({\wt \Pi}_c^{n}c^{n}-c^{n})\|_{L^2})  
  \|{\wt \theta}_c^{n+1}\|_{L^6} 
  \\
  \le &C (\tau+h)
  (\|\nabla {\wt \theta}_c^{n}\|_{L^2}+h^{r})  
  \|{\wt \theta}_c^{n+1}\|_{L^6} 
  \\
  \le & (C\tau+Ch+\epsilon)(\|{\wt \theta}_c^{n}\|_{H^1}^2  
  +\|{\wt \theta}_c^{n+1}\|_{H^1}^2)+C\epsilon^{-1}(h^{2r+2}+\tau^2),
\end{align*} 
and
\begin{align*} 
  & \big |\big(( \U_h^{n+1}- \u^{n+1})\cdot\nabla c^{n},\,  
  {\wt \theta}_c^{n+1}\big)\big | 
  \\
  \le& |\big({\wt \theta}_\u^{n+1}\cdot\nabla c^{n},\,  
  {\wt \theta}_c^{n+1}\big)\big |+
  |\big(({\wt \u}_h^{n+1}- \u^{n+1})\cdot\nabla c^{n},\,  
  {\wt \theta}_c^{n+1}\big)\big |
  \\
  \le& C\|{\wt \theta}_u^{n+1}\|_{L^2}\|{\wt \theta}_c^{n+1}\|_{L^2}
  +C\|{\wt \u}_h^{n+1}- \u^{n+1}\|_{H^{-1}} 
  \|\nabla c^{n}{\wt \theta}_c^{n+1}\|_{H^1}
  \\ 
  \le& \epsilon  
  \|{\wt \theta}_c^{n+1}\|_{H^1}^2+C\|{\wt \theta}^{n}_c \|_{L^2}^2
  +C\|{\wt \theta}_c^{n+1} \|_{L^2}^2+C\epsilon^{-1}(h^{2r+2}+\tau^2).
\end{align*} 
It follows that 
\begin{align*}
  |J_3({\wt \theta}_c^{n+1})|
  \le 
  (C\tau+Ch +\epsilon )(\|{\wt \theta}_c^{n}\|_{H^1}^2+\|{\wt \theta}_c^{n+1}\|_{H^1}^2 )
  +C\|{\wt \theta}^{n}_c \|_{L^2}^2+C\|{\wt \theta}_c^{n+1} \|_{L^2}^2+C\epsilon^{-1}(h^{2r+2}+\tau^2). 
\end{align*}
Therefore, 
\begin{align}
  &\frac{1}{2} D_t \|{\wt \theta}_c^{n+1}\|^2_{L^2} +
  \big \| \sqrt{D(\U_h^{n+1})}\nabla{\wt \theta}_c^{n+1}\big \|^2_{L^2}\nn\\
  \leq&
  C (\epsilon+ \tau +h)
  ( \|\nabla {\wt \theta}_c^{n}\|_{L^2}^2 + \|\nabla {\wt \theta}_c^{n+1}\|_{L^2}^2 ) 
  + C\epsilon^{-1}(\|{\wt \theta}_c^{n+1}\|_{L^2}^2+\|{\wt \theta}^{n}_c\|_{L^2}^2)
  + C \epsilon^{-1}(h^{2r+2}+\tau^2)\nn\\
  &~~~+ C\big \|D_t (c^{n+1}-{\wt \Pi}_c^{n+1} c^{n+1})
  \big \|_{L^2}^2 + C\|T_c^{n+1}\|_{L^2}^2  \, .
  \label{ntem-1}
\end{align}
By Gronwall's inequality and \refe{np2-3},
there exists $\tau_0>0$ and $h_0>0$ such that,
\begin{align}
  \|{\wt \theta}_c^{n+1}\|_{L^2} \leq C(\tau + h^{r+1})
  \label{ngw3-1}
\end{align}
for $0\leq n\leq k$, when $\tau<\tau_0 \le \widehat \tau_0$ 
and $h<h_0\le \widehat h_0$.

Similarly to the estimate for $\theta_p^{n}$ in section 3.1, 
we can get the estimate to ${\wt \theta}_p^{n}$
\begin{align}
  \|{\wt \theta}_p^{n+1}\|_{L^2} \leq C(\tau + h^{r+1}).
  \label{nPp-Ph}
\end{align}
\refe{nerror-1} follows \refe{nUp-Uh}, \refe{ngw3-1} and \refe{nPp-Ph}. 
Theorem \ref{ErrestFEMSol} is proved by combining 
\refe{nU-Up-Lp}-\refe{nU-Up-H-1}, \refe{nerror-1} and the basic 
projection error estimates in Lemma \ref{nEpro}.  
The proof is complete. \quad \endproof


\section{Proof to Lemma 3.2}
\setcounter{equation}{0}
To prove Lemma \ref{EPro}, we define an extra classical elliptic projection 
$R_h^{n}$: $H^1(\Omega) \rightarrow V_h^{r}$ by
\begin{align}
  &\Big(D(\U^{n})\nabla (R_h^{n}\psi-\psi), \nabla \phi_h \Big)  = 0 ,
  \qquad \qquad \mbox{for~all}~~\phi_h\in V_h^{r},\label{Pro2}
\end{align}
with $\int_\Omega(\psi-R_h^{n}\psi)d x=0$. 
Similarly, by classical FE theory \cite{BS,VThomee},
\begin{align}
  \|R_h^{n}\psi\|_{W^{1,p}}\le C \|\psi\|_{W^{1,p}}, \quad 1 < p \le 4. 
  \label{ph}
\end{align}

For $2\leq p \leq 4$ and $1/q+1/p=1$, 
by \refe{Pro2} we have 
\begin{align}
  & |(D(\U^{n})\nabla(\Pi_c^{n}{\cal C}^{n}-{\cal C}^{n}),\nabla v)| 
  \nn \\ 
  \le &
  |(D(\U^{n})\nabla(\Pi_c^{n}{\cal C}^{n}-I_h{\cal C}^{n}), 
  \nabla v)|
  + |(D(\U^{n})\nabla(I_h{\cal C}^{n}-{\cal C}^{n}),\nabla v)| 
  \nn \\
  \le& 
  |(D(\U^{n})\nabla(\Pi_c^{n}{\cal C}^{n}-I_h {\cal C}^{n}), 
  \nabla R_h^{n}v)| 
  +C\|\nabla(I_h{\cal C}^{n}-{\cal C}^{n})\|_{L^p} 
  \|\nabla v\|_{L^q} 
  \nn \\
  \le & 
  |(D(\U^{n})\nabla(\Pi_c^{n}{\cal C}^{n}-{\cal C}^{n}), 
  \nabla R_h^{n} v)|
  + C h (\|\nabla v\|_{L^q} + \|\nabla R_h^{n} v\|_{L^q}),
  \nn 
\end{align} 
where we have noted the fact 
that $D(\U^n)$ and $D({\wt\U}_h^n)$ are symmetric.
With \refe{Nonclassical} and \refe{ph} we further have 
\begin{align} 
  & |(D(\U^n)\nabla(\Pi_c^{n}{\cal C}^{n}-{\cal C}^{n}),\nabla v)| 
  \nn \\ 
  \le&  
  |((D(\U^n)-D({\wt\U}_h^n)) \nabla\Pi_c^{n+1}{\cal C}^{n+1}, 
  \nabla R_h^{n}v )|
  + Ch \|\nabla v\|_{L^q}  
  \nn \\
  \le & 
  C \left ( 
    \|\U^n-{\wt\U}_h^n \|_{L^p} 
    \| \nabla {\cal C}^{n} \|_{L^\infty} 
    + 
    \|\U^n-{\wt\U}_h^n \|_{L^\infty} 
  \| \nabla (\Pi_c^{n}{\cal C}^{n} - {\cal C}^{n}) \|_{L^p} + h \right ) 
  \|\nabla v\|_{L^q}  
  \nn \\
  \le&  
  C \left ( h + h^{1/4} 
  \| \nabla (\Pi_c^{n+1}{\cal C}^{n+1} - {\cal C}^{n+1}) \|_{L^p} \right ) 
  \|\nabla v\|_{L^q} ,
  \label{w1p}
\end{align}
where we have also used some basic estimates to $R_h^{n}$. 
It follows that 
\begin{align} 
  \|\nabla(\Pi_c^{n}{\cal C}^{n}-{\cal C}^{n})\|_{L^{p}}
  \le 
  C \|D(\U^{n})\nabla(\Pi_c^{n}{\cal C}^{n}-{\cal C}^{n})\|_{L^{p}} 
  \le C h,
  \label{Pi2-p} 
\end{align} 
when $h \le h_1$ for some $h_1>0$. By using an inverse inequality, we see that 
\begin{align}
  \|\Pi_c^{n}{\cal C}^{n}\|_{W^{1,\infty}} 
  \leq&\|I_h{\cal C}^{n}\|_{W^{1,\infty}} 
  +\|\Pi_c^{n}{\cal C}^{n}- I_h {\cal C}^{n}\|_{W^{1,\infty}}\nn \\
  \leq&Ch^{-d/4}\|\Pi_c^{n}{\cal C}^{n}-I_h{\cal C}^{n}\|_{W^{1,4}} 
  +C\nn \\
  \leq&Ch^{1-d/4}+C.\nn 
\end{align}
We have proved \refe{p2-2} and the second part of \refe{p2-1}. 

To get the $L^2$-norm estimate in \refe{p2-1}, we use the duality argument and 
consider the equation 
\begin{align}
  \label{bc-w}
  \left.
    \begin{aligned}
      &-\nabla\cdot\Big(D(\U^{n})\nabla w\Big)  
      ={\cal C}^{n}-\Pi_c^{n}{\cal C}^{n}, 
      \qquad &&\mbox{in } \Omega, 
      \\
      &\hskip1in D(\U^{n})\nabla w\cdot {\bf n} = 0,  
      \qquad &&\mbox{on } \partial\Omega,
    \end{aligned}
  \right.
\end{align} 
with $\int_{\Omega}wdx = 0$. 
Its solution satisfies 
\begin{align} 
  \| w \|_{H^2}  
  \le C \| {\cal C}^{n}-\Pi_c^{n}{\cal C}^{n} \|_{L^2} \, . 
  \label{h2} 
\end{align} 
By \refe{Nonclassical}, 
\begin{align}
  &\|{\cal C}^{n}-\Pi_c^{n}{\cal C}^{n}\|_{L^2}^2 
  \nn 
  \\ 
  =& \Big(D(\U^n)\nabla (w-I_h w), 
  \nabla ({\cal C}^{n}-\Pi_c^{n}{\cal C}^{n})\Big)
  \nn \\ 
  & +\Big(\nabla (I_hw), 
  D(\U^n)\nabla ({\cal C}^n-\Pi_c^{n}{\cal C}^{n})\Big) 
  \nn \\
  = &\Big(D(\U^n)\nabla (w-I_h w), 
  \nabla ({\cal C}^{n}-\Pi_c^{n}{\cal C}^{n})\Big)
  \nn \\ 
  &+ \Big((D({\wt\U}_h^{n})-D(\U^{n}))\nabla \Pi_c^{n}{\cal C}^{n}, 
  \nabla I_h w \Big) 
  \nn \\
  = &\Big(D(\U^n)\nabla (w-I_h w), 
  \nabla ({\cal C}^{n}-\Pi_c^{n}{\cal C}^{n})\Big)
  \nn \\ 
  &+ \Big((D({\wt\U}_h^{n})-D(\U^{n}))\nabla {\cal C}^{n}, 
  \nabla I_h w \Big) 
  \nn \\
  &+ \Big((D({\wt\U}_h^{n})-D(\U^{n}))(\nabla \Pi_c^{n}{\cal C}^{n}-\nabla{\cal C}^{n}), 
  \nabla I_h w \Big) 
  \label{p3-l2} 
\end{align} 
For the second term in the right hand side of the last equation, we have 
the following estimate 
\begin{align} 
  & |((D({\wt\U}_h^n)-D(\U^n))\nabla {\cal C}^{n}, 
  \nabla I_h w ) |
  \nn \\ 
  \le&   
  |((D({\wt\U}_h^n)-D(\U^n))\nabla {\cal C}^{n}, \nabla w)|  
  + 
  |((D({\wt\U}_h^{n})-D(\U^n))\nabla {\cal C}^{n},  
  \nabla (I_hw-w) )| \nn \\
  \le& C\|{\wt\U}_h^{n}-\U^{n}\|_{H^{-1}}\|\nabla {\cal C}^{n}\nabla w\|_{H^1} 
  + C\|{\wt\U}_h^n-\U^n\|_{L^{2}}\|\nabla (I_hw-w)\|_{L^2} \| \nabla {\cal C}^n \|_{L^\infty} \nn \\
  \le& Ch^2\|w\|_{H^2}.
  \label{4.6-1}
\end{align} 
where we have used \refe{U-Up-Lp}-\refe{U-Up-H-1}. 
\begin{align}
  \|{\cal C}^{n}-\Pi_c^{n}{\cal C}^{n}\|_{L^2}^2 
  \le& C h\|w \|_{H^2}  
  \| \nabla ({\cal C}^{n}-\Pi_c^{n}{\cal C}^{n}) \|_{L^2} 
  + Ch^2\|{\cal C}^{n}-\Pi_c^{n}{\cal C}^{n}\|_{L^2}
  \nn \\
  & + \|{\wt\U}_h^{n}-\U^{n}\|_{L^3} 
  \|\nabla(\Pi_c^{n}{\cal C}^{n}-{\cal C}^{n})\|_{L^2} 
  \|\nabla I_h w\|_{L^6}
  \nn \\ 
  \le &
  Ch^2\|{\cal C}^{n}-\Pi_c^{n}{\cal C}^{n}\|_{L^2}
  \, , 
  \nn 
\end{align} 
and \refe{p2-1} follows immediately. 

It remains to prove \refe{p2-3}. From \refe{Pro2}, we see that 
\begin{align}
  &\Big(D(\U^n)\nabla D_{\tau}(\Pi_c^{n}{\cal C}^{n} 
  -{\cal C}^{n}), \nabla \phi_h \Big) 
  +\Big(D_{\tau}(D(\U^{n}))\nabla (\Pi_c^{n-1}{\cal C}^{n-1}-{\cal C}^{n-1}),  
  \nabla \phi_h \Big)  
  \nn \\
  &+\Big((D({\wt\U}_h^{n})-D(\U^n))\nabla(D_{\tau}\Pi_c^{n}{\cal C}^{n}),  
  \nabla \phi_h \Big) 
  +\Big(D_{\tau}(D({\wt\U}_h^{n})-D(\U^n))\Pi_c^{n-1}\nabla{\cal C}^{n-1}, 
  \nabla \phi_h \Big) = 0 .\nn 
\end{align}
Using similar proof for 
\refe{p2-1}, we can get the desired result.
The proof is complete.
\quad 
\endproof 
\section{Numerical examples}
\setcounter{equation}{0}
In this section, we present numerical results
for incompressible miscible flows in both two and three-dimensional 
porous media to confirm our theoretical analysis and show the efficiency 
of linearized Galerkin FEMs. 
Here we always assume that the solution of the system is smooth. 
The problem with non-smooth solutions was considered in \cite{BJM}, where
the convergence of a discontinuous-mixed FEM was proved. 
All computations in this section 
are performed by using the software FEniCS \cite{fenics}.

\medskip

We rewrite the system (\ref{e-1})-(\ref{e-4}) by
\begin{align}
  &\frac{\partial c}{\partial t}-\nabla\cdot(D(\u)\nabla c)
  +\u\cdot\nabla c= g,
  \label{n-e-1}\\[3pt]
  &\nabla\cdot\u=f,
  \label{n-e-2}\\[3pt]
  &\u=-\frac{1}{\mu(c)}\nabla p,
  \label{n-e-3}
\end{align}
where $D(\u)=1+|\u|^2/(1+|\u|^2)+\u\otimes\u$ and $\mu(c)=1+c^2$.
\vskip0.1in 

{\bf Example 5.1.}
First, we consider a two-dimensional model where $\Omega = [0,1] \times [0,1]$.
We set the terminal time $T = 1.0$. The functions $f$ and $g$
are chosen correspondingly to the exact solution
\begin{align}
  &p=1+1000x^2(1-x)^3y^2(1-y)^3t^2e^{-t},
  \label{num001}\\
  &c=0.2+50x^2(1-x)^2y^2(1-y)^2te^{t},
  \label{num002}
\end{align}
which satisfies the boundary condition (\ref{e-3}).
\begin{figure}[ht]
  \centering
  \includegraphics[width = 65mm]{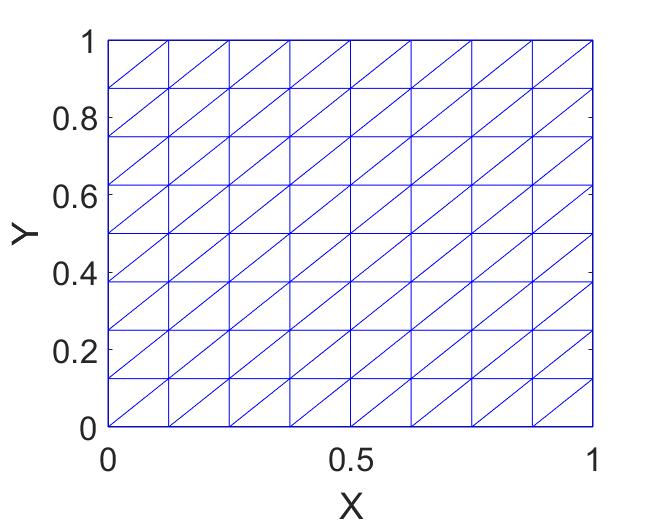}
  \caption{A uniform triangular mesh on the unit square domain with M = 8}
  \label{mesh}
\end{figure}

\begin{table}[h]
  \centering
  \begin{center}
    \caption{$L^2$-norm errors of linearized Galerkin-mixed 
      FEM \refe{e-FEM-1}-\refe{e-FEM-3} in 2D 
    (Example 5.1). }
    \vskip 0.1in
    \label{linear-L2-1}
    \begin{tabular}{c|ccc}
      \hline 
      $\tau = \frac{8}{M^2}(r=1)$&$Err_P$ &$Err_\U$&$Err_{\cal C}$\\\hline
      M = 8        	              &2.63e-02&1.99e-01     &5.09e-02     \\
      M = 16  	                  &1.29e-02&1.01e-01     &1.20e-02     \\ 
      M = 32  	                  &6.38e-03&5.07e-02     &2.93e-03     \\ 
      M = 64  	                  &3.18e-03&2.54e-02     &7.29e-04     \\
      M =128  	                  &1.59e-03&1.27e-02     &1.82e-04     \\\hline
      order 	                  &1.01    &0.99         &2.03         \\\hline
      \hline 
      $\tau = \frac{16}{M^3}(r=2)$&$Err_P$ &$Err_\U$&$Err_{\cal C}$\\\hline
      M = 8        	              &3.48e-03&2.81e-02     &4.66e-03     \\
      M = 16  	                  &8.53e-04&7.23e-03     &5.64e-04     \\ 
      M = 32  	                  &2.12e-04&1.82e-03     &6.94e-05     \\
      M = 64  	                  &5.30e-05&4.56e-04     &8.62e-06     \\
      M =128  	                  &1.33e-05&1.14e-04     &1.08e-06     \\\hline
      order 	                  &2.01    &1.99         &3.02         \\\hline
      \hline 
    \end{tabular}
  \end{center}
\end{table}

We use a uniform triangular mesh with M+1 vertices in each direction, where 
$h = \frac{\sqrt{2}}{M}$ (see Figure \ref{mesh} for  
the illustration with $M=8$). 
We solve the system \refe{n-e-1}-\refe{n-e-3}
by the scheme \refe{e-FEM-1}-\refe{e-FEM-3} with $r=1$ and $r=2$, respectively.  

As the expected optimal convergence rate in $L^2$-norm 
is $O(\tau+h^{r+1})$, we set  $\tau = \frac{8}{M^{r+1}}$ 
in our computation. 
The $L^2$-norm errors of $P_h^{N}$, $\U_h^{N}$ and ${\cal C}_h^{N}$ are presented in 
Table \ref{linear-L2-1} for $r=1,2$,
where $Err_v = \|v_h^{N}-v(\cdot,t_N)\|_{L^2}$. 
From Table \ref{linear-L2-1},
we can see clearly that the convergence rates for $P_h^{N}$ 
and $\U_h^{N}$ are optimal with the order $O(h^{r})$, 
while the rate for ${\cal C}_h^{n}$ is optimal with the order $O(h^{r+1})$ 
for both $r=1,2$. 
The one-order lower approximation to $(p, \u)$ 
does not affect the accuracy of the numerical concentration and the mesh-size restriction 
in \refe{mesh-cond-2} is not necessary.


\vskip0.1in 

{\bf Example 5.2.}
Secondly, we study the system \refe{n-e-1}-\refe{n-e-3}
in a three-dimensional cube $[0,1] \times [0,1] \times [0,1]$., where 
the functions $f$ and $g$ are chosen correspondingly to the smooth 
exact solution
\begin{align}
  &p=1.0+1000x^2(1-x)^3y^2(1-y)^3z^2(1-z)^3t^2e^{-t},
  \label{num003}\\
  &c=0.2+50x^2(1-x)^2y^2(1-y)^2z^2(1-z)^2te^{t}.
  \label{num004}
\end{align}
We also set the terminal time T = 1.0 in this example.

A uniform triangular mesh with M+1 vertices in each direction 
of the cube is used in our computation, 
where $h = \frac{\sqrt{3}}{M}$. We solve the system  
by the linearized Galerkin FEMs in 
\refe{e-FEM-1}-\refe{e-FEM-3} 
by the lowest-order Galerkin-mixed FEM ($r=1$). Such a Galerkin-mixed FEM 
is most popular in practical computation, 
particular for problems with discontinuous media. 
To show the accuracy in spatial direction, 
we set $\tau = \frac{8}{M^2}$ in our computation. 
We present in Table \ref{table3} 
the $L^2$-norm errors of concentration, velocity and pressure. 
Again numerical results confirm our theoretical analysis and 
the accuracy of the numerical concentration is in the order $O(h^2)$, while 
all previous analyses only showed the first-order accuracy of numerical 
concentration. 

\begin{table}[!ht]
  \caption{$L^2$-norm errors of linearized Galerkin FEM \refe{e-FEM-1}-\refe{e-FEM-3} 
  in 3D (Example 5.2)}
  \vskip0.1in 
  \centering
  \begin{tabular}{c|ccc}
    \hline 
    $\tau = \frac{8}{M^2}$&$Err_P$ &$Err_\U$&$Err_{\cal C}$\\\hline
    M = 8        	       &5.70e-04&5.36e-03     &9.05e-04      \\
    M = 16  	           &2.82e-04&2.72e-03     &2.40e-04      \\
    M = 32  	           &1.40e-04&1.36e-03     &6.10e-05      \\
    M = 64  	           &7.55e-05&7.13e-04     &1.38e-05      \\\hline
    order 	               &0.98    &0.97         &2.01          \\\hline
    \hline 
  \end{tabular}
  \label{table3}
\end{table}

\section{Conclusions}
We have presented unconditionally optimal error analysis of 
commonly-used Galerkin-mixed FEMs for a nonlinear and strongly 
coupled parabolic system from incompressible miscible flow in 
porous media. In particular, for the most popular lowest-order 
Galerkin-mixed method, we show unconditionally the second-order 
accuracy $O(h^2)$ for the numerical concentration ${\cal C}_h^n$ 
in spatial direction. In all previous works, only the first-order 
accuracy was obtained under certain time-step restriction and the 
mesh-size condition. Moreover, our analysis is based on a 
quasi-projection and the approach presented in this paper can be 
extended to many other coupled nonlinear parabolic PDEs to obtain 
optimal error estimates for all components. With the numerical 
concentration $C_h^n \in V_h^r$ obtained by the lowest-order 
Galerkin-mixed FEM, one can get new numerical velocity and pressure 
of the accuracy $O(h^2)$ (same as the concentration) by resolving 
the elliptic pressure equation \refe{e-3} at a given time with a 
higher-order Galerkin FEM or a higher-order mixed method. 

{
  Optimal error estimates presented in this paper is based on strong
  regularity assumptions of the solution of the system and physical 
  parameters as usual. Since the present paper focuses on the new 
  analysis of classical Galerkin-mixed FEMs, the regularity 
  assumptions in $(A1)-(A6)$ may not be optimal. Some related works 
  under weaker assumptions were done by several authors \cite{CLLS,LS3,RKN}. 
  Also the existence and 
  uniqueness of the strong solution of the system 
  \refe{e-1}-\refe{e-2} in a general setting remains open. A 
  systematic numerical simulation on the above scheme with many 
  other approximations will be presented in the coming article 
  \cite{WS} and the problems with non-smooth domains and 
  discontinuous coefficients will be also studied there. 

}
\bigskip

\noindent{\bf Acknowledgements}~The authors would like to thank the 
anonymous referees for their valuable suggestions and comments.


\bigskip

\begin{thebibliography}{99}
  \bibitem{AFS}
    {
      A. Agosti, L. Formaggia and A. Scotti,
      Analysis of a model for precipitation and dissolution coupled 
      with a Darcy flux,
      {\em J. Math. Anal. Appl.}, 431(2015), 752--781.
    }
  \bibitem{AO}
    B.~Amaziane and M.~El Ossmani,  
    Convergence analysis of an approximation to miscible fluid flows 
    in porous media by combining mixed finite element and finite volume methods, 
    {\em Numer. Methods Partial Diff. Eq.}, 24(2008), 799--832.

  \bibitem{AS}
    R. An and J. Su, 
    Optimal error estimates of semi-implicit Galerkin method for time-dependent 
    Nematic liquid crystal flows, 
    {\em J. Scientific Computing}, 74(2018), 979--1008. 

  \bibitem{AW}
    {
      T. Arbogast and W. Wang,
      Stability, monotonicity, maximum and minimum principles, and 
      implementation of the volume corrected characteristic method,
      {\em SIAM J. Sci. Comput.}, 33(2011), 1549--1573.
    }

  \bibitem{AWZ}
    {
      T. Arbogast, M.F. Wheeler and N. Zhang,
      A nonlinear mixed finite element method for a degenerate parabolic
      equation arising in flow in porous media,
    {\em SIAM J. Numer. Anal.}, 33(1996), 1669--1687.}

  \bibitem{BK}
    {
      J.W. Barrett and P. Knabner,
      Finite element approximation of the transport of reactive solutes
      in porous media. Part 1: error estimates for nonequilibrium
      adsorption processes,
      {\em SIAM J. Numer. Anal.}, 34(1997), 201--227.
    }
  \bibitem{BJM}
    S.~Bartels, M.~Jensen and R.~Muller,
    Discontinuous Galerkin finite element convergence for 
    incompressible miscible displacement problems
    of low regularity,
    {\em SIAM J. Numer. Anal.}, 47(2009), 3720--3742.

  \bibitem{BB}
    J.~Bear and Y.~Bachmat, 
    Introduction to Modeling of Transport Phenomena in Porous Media, 
    Springer-Verlag, New York, 1990.


  \bibitem{BS} 
    S.~Brenner and L.~Scott,
    The Mathematical Theory of Finite Element Methods,
    Springer, New York, 2002.

  \bibitem{CGW}
    W. Cai, Max Gunzburger and J. Wang, 
    Convergence analysis of Crank-Nicolson
    Galerkin-Galerkin FEMs for miscible displacement
    in porous media, to appear. 

  \bibitem{CLLS}
    {W. Cai, B. Li, Y. Lin and W. Sun,
      Analysis of fully discrete FEM for miscible displacement in porous
      media with Bear--Scheidegger diffusion tensor,
    {\em Numer. Math.}, 141(2019), 1009--1042.}

  \bibitem{CCW}
    F.~Chen, H.~Chen and H.~Wang, 
    An optimal-order error estimate for a Galerkin-mixed  
    finite-element time-stepping procedure for porous media flows, 
    {\em Numer. Methods Partial Differ. Equations}, 28.2(2012), 707--719. 

  \bibitem{CE}
    {Z. Chen and R. Ewing, 
      Mathematical analysis for reservoir models, 
      {\em SIAM J. Math. Anal.}, 30 (1999), 43--453. 
    } 

  \bibitem{CHM}
    {Z. Chen, G. Huan and Y. Ma, 
      Computational Methods for Multiphase Flows in Porous Media, 
      Computational science and engineering,  SIAM, PA, 2006.
    }   

  \bibitem{CWW}
    A.~Cheng, K.~Wang and H.~Wang, 
    Superconvergence for a time-discretization procedure 
    for the mixed finite element
    approximation of miscible displacement in porous media,
    {\em Numer. Methods Partial Differ. Equations}, 28(2012), 1382--1398.

  \bibitem{CH}
    {K. Chrysafinos and L.S. Hou, 
      Error estimates for semidiscrete finite element approximations 
      of linear and semilinear parabolic equations under minimal 
      regularity assumptions, 
    {\em SIAM J. Numer. Anal.}, 40(2002), 282--306.}


  \bibitem{DEW1}
    J. Douglas, JR., R.E.~Ewing and M.F.~Wheeler, 
    The approximation of the pressure by a
    mixed method in the simulation of miscible displacement,  
    {\em RAIRO Analyse num{\'e}rique}, 17(1983), 17--33.


  \bibitem{DEW2}
    J. Douglas, JR., R. Ewing and M.F. Wheeler, 
    A time-discretization procedure for a mixed finite element 
    approximation of miscible displacement in porous media, 
    {\em RAIRO Anal. Numer.}, 17(1983), 249--265.

  \bibitem{DFP}
    J. Douglas, JR., F. Furtada and F. Pereira, 
    On the numerical simulation of waterflooding of 
    heterogeneous petroleum reservoirs,
    {\em Comput. Geosciences}, 1(1997), 155--190.

  \bibitem{JJ}
    J. Douglas, JR. and J.E. Roberts,  
    Global estimates for mixed methods for second order 
    elliptic equations,
    {\em Math. Comput.}, 44(1985), 39--52.

  \bibitem{Dur}
    R.G. Dur\'{a}n, 
    On the approximation of miscible displacement in
    porous media by a method of characteristics combined with a mixed
    method, 
    {\em SIAM J. Numer. Anal.}, 25(1988), 989--1001.

  \bibitem{DR}
    R.G. Dur\'{a}n, 
    Error analysis in $L^p$, $1\leq p \leq \infty$, for 
    mixed finite element methods for linear and quasi-linear 
    elliptic problems, 
    {\em RAIRO Mod. Math. Anal. Num\'{e}r.}, 22(1988), 371--387.

  \bibitem{EM}
    V.J.~Ervin and W.W.~Miles, 
    Approximation of time-dependent viscoelastic fluid flow: 
    SUPG approximation, 
    {\em SIAM J. Numer. Anal.}, 41(2003), 457--486.

  \bibitem{ER}
    {R.E.~Ewing and T.F.~Russell,
      Efficient time-stepping methods for miscible displacement problems
      in porous media,
    {\em SIAM J. Numer. Anal.}, 19(1982), 1--67.}

  \bibitem{EW}
    R.E.~Ewing and M.F.~Wheeler, 
    Galerkin methods for miscible displacement problems in porous media, 
    {\em SIAM J. Numer. Anal.}, 17(1980), 351--365.

  \bibitem{Ewing}
    R.E. Ewing, ed,
    The mathematics of Reservoir Simulation,
    Frontiers in Applied Mathematics, SIAM, Philadelphia, PA, 1983.

  \bibitem{ERW2}
    R.E.~Ewing, T.F.~Russell and M.F.~Wheeler, 
    Convergence analysis of an approximation of miscible displacement 
    in porous media by mixed finite elements and a modified method 
    of characteristics,
    {\em Comput. Methods Appl. Mech. Engrg.}, 47(1984), 73--92.

  \bibitem{EWang}
    R.E. Ewing and H. Wang,
    A summary of numerical methods for time-dependent advection-dominated
    partial differential equations,
    {\em J. Comput. Appl. Math.}, 128(2001), 423--445.

  \bibitem{Feng}
    X.~Feng, 
    On existence and uniqueness results for a coupled system
    modeling miscible displacement in porous media, 
    {\em J. Math. Anal. Appl.}, 194(1995), 883--910.


  \bibitem{Feng2}
    {X.~Feng, 
      Recent developments on modeling and analysis of flow of miscible
      fluids in porous media, 
    {\em Contemp. Math.}, 295(2002), 229--240.}

  \bibitem{FN}
    X.~Feng and M.~Neilan,
    A modified characteristic finite element method for a fully nonlinear
    formulation of the semigeostrophic flow equations,
    {\em SIAM J. Numer. Anal.},  47(2009), 2952--2981.

  \bibitem{KNP}
    S.~Kumar, N.~Nataraj and A.K.~Pani,
    Finite volume element method for the incompressible miscible displacement
    problems in porous media,
    {\em Proc. Appl. Math. Mech.}, 7(2007), 2020015--2020016.

  \bibitem{KY}
    {
      S. Kumar and S. Yadav,
      Modified method of characteristics combined with finite volume
      element methods for incompressible miscible displacement problems
      in porous media.
      {\em Int. J. PDEs}, 2014(2014). 
    }


  \bibitem{LS2}
    B.~Li and W. Sun,
    Unconditional convergence and optimal error estimates of a
    Galerkin-mixed FEM for incompressible miscible flow in porous media,
    {\em SIAM J. Numer. Anal.}, 51(2013), 1959--1977.

  \bibitem{LS3} 
    B. Li and W. Sun,
    Regularity of the diffusion-dispersion
    tensor and error analysis of FEMs for a porous media flow,
    {\em SIAM J. Numer. Anal.}, 53(2015), 1418--1437.


  \bibitem{fenics}
    A.~Logg, K.~Mardal and G.~Wells (Eds.),
    Automated Solution of Differential Equations by the
    Finite Element Method,   
    Springer, Berlin, 2012.

  \bibitem{MLG}
    S.M.C.~Malta, A.F.D.~Loula and E.L.M~Garcia, 
    Numerical analysis of a stabilized finite element method 
    for tracer injection simulations, 
    {\em Comput. Methods Appl. Mech. Engrg.}, 187(2000), 119--136. 

  \bibitem{RKN}
    {
      F.A. Radu, K. Kumar, J.M. Nordbotten {and I.S. Pop},
      A robust, mass conservative scheme for two-phase flow
      in porous media including H\"older continuous nonlinearities,
      {\em IMA J. Numer. Anal.}, 38(2018), 884--920.
    }
  \bibitem{RPA}
    {
      F.A. Radu, I.S. Pop and S. Attinger,
      Analysis of an Euler implicit-mixed finite element scheme for 
      reactive solute transport in porous media,
      {\em Numer. Methods Partial Differential Equations}, 26(2010), 320--344.
    }

  \bibitem{RT} 
    P.A. Raviart and J.M. Thomas, 
    A mixed finite element method
    for 2nd order elliptic problems, Mathematical Aspects of
    Finite Element Methods, Lecture Notes in Math., vol. 606,
    Springer-Verlag, 1977, 292--315.

  \bibitem{RW}
    {
      B.M. Rivi\`ere and N.J. Walkington,
      Convergence of a discontinuous Galerkin method for the 
      miscible displacement equation under low regularity,
      {\em SIAM J. Numer. Anal.}, 49(2011), 1085--1110.
    }

  \bibitem{Rus1} 
    T. F. Russell,
    Time stepping along characteristics with incomplete iteration for a Galerkin
    approximation of miscible displacement in porous media, 
    {\em SIAM J. Numer. Anal.}, 22(1985), 970--1013.


  \bibitem{SWML}
    G.~Scovazzi, M.F.~Wheeler, A.~Mikeli\'{c} and S.~Lee,
    Analytical and variational numerical methods for unstable miscible 
    displacement flows in porous media,
    {\em J. Comput. Phys.}, 335(2017), 444--496.

  \bibitem{VThomee}
    V. Thom\'{e}e, 
    Galerkin finite element methods for parabolic problems, 
    Springer-Verkag Berkub Geudekberg 1997.

  \bibitem{Wang}
    H.~Wang, 
    An optimal-order error estimate for a family of
    ELLAM-MFEM approximations to porous medium flow, 
    {\em SIAM J.  Numer. Anal.}, 46(2008), 2133--2152.

  \bibitem{WLELQ}
    H. Wang, D. Liang, R.E. Ewing, S.L. Lyons and G. Qin,  
    An approximation to miscible fluid flows in porous media with point sources 
    and sinks by an Eulerian-Lagrangian localized adjoint method and mixed 
    finite element methods, 
    {\em SIAM J. Sci. Comput.}, 22(2000), 561--581. 

  \bibitem{WSS}
    J. Wang, Z. Si and W. Sun,
    A new error analysis of characteristics-mixed FEMs for miscible
    displacement in porous media,
    {\em SIAM J. Numer. Anal.}, 52(2014), 3300--3020.

  \bibitem{Whe}
    M.F. Wheeler, 
    A priori $L^2$ error estimates for Galerkin
    approximations to parabolic partial differential equations,
    {\em SIAM J. Numer. Anal.}, 10(1973), 723--759.

  \bibitem{WS}
    C.~Wu and W.~Sun, 
    Efficient fully-discrete Galerkin-mixed FEMs for incompressible 
    miscible flow in porous media, in preparation. 



  \bibitem{ZYS}
    H.~Zheng, J.~Yu and L.~Shan,  
    Unconditional error estimates for time dependent viscoelastic fluid flow, 
    {\em Appl. Numer. Math.}, 119(2017), 1--17. 

\end{thebibliography}
\end{document}